\documentclass[reqno,11pt,noamsfonts]{amsart}
\usepackage[foot]{amsaddr}         

\usepackage[default,scale=0.92]{opensans}
\usepackage[notext]{stix}

\usepackage{algorithm}               
\usepackage{algpseudocode}  

\usepackage{etoolbox}
\patchcmd{\section}{\scshape}{\bfseries}{}{}
\makeatletter
\renewcommand{\@secnumfont}{\bfseries}
\makeatother

\usepackage{amsthm}
\newtheorem{remark}{Remark}

\usepackage{bm}
\usepackage{mathtools}
\usepackage{xfrac}                  
\usepackage{accents}                

\usepackage[svgnames]{xcolor}

\usepackage{graphicx}
\usepackage[labelfont=bf]{caption}
\usepackage{subcaption}
\usepackage{hyperref}
\usepackage{enumitem}
\usepackage{array}             
\usepackage{booktabs}          

\usepackage[numbers,sort&compress]{natbib}

\calclayout

\begin{document}




\newcommand   {\mathsc}[1]{\text{{\rmfamily\scshape #1}}}
\renewcommand {\mathit}[1]{\text{{\rmfamily\itshape #1}}}

\newcommand {\CFL}    {\mathit{CFL}}

\newcommand {\EX}     {\mathrm{ex}}       
\newcommand {\IM}     {\mathrm{im}}       

\newcommand{\vs} {\nu_{\mathrm s}}
\newcommand{\cs} {c_{\mathrm s}}
\newcommand{\ds} {\kappa_{\mathrm s}}

\newcommand {\fc}     {\M f_\mathrm{c}}
\newcommand {\hc}     {\M h_\mathrm{c}}
\newcommand {\fd}     {\M f_\mathrm{d}}
\newcommand {\fs}     {\M f_\mathrm{s}}
\newcommand {\fex}    {\M f_\mathrm{ex}}
\newcommand {\fim}    {\M f_\mathrm{im}}

\newcommand {\Ac}     {\M A_\mathrm{c}}
\newcommand {\Ad}     {\M A_{\mathrm{d}}}
\newcommand {\nud}    {\tilde\nu}
\newcommand {\kappad} {\tilde\kappa}

\newcommand {\PO}  {\mathcal M}
\newcommand {\g}   {\mathcal G}
\newcommand {\F}   {\mathcal F}
\newcommand {\Fc}  {\mathcal F_{\!\mathrm{c}}}
\newcommand {\Fd}  {\mathcal F_{\!\mathrm{d}}}
\newcommand {\Fs}  {\mathcal F_{\!\mathrm{s}}}
\newcommand {\Fex} {\mathcal F_{\!\mathrm{ex}}}
\newcommand {\Fim} {\mathcal F_{\!\mathrm{im}}}

\newcommand {\lbf} {{\ell\!}}

\newcommand {\iu}  {\mathrm{i}}
\newcommand {\zi}  {z_\mathrm{i}}
\newcommand {\zr}  {z_\mathrm{r}}

\newcommand {\V}  [1] {\vec{#1}}                       
\newcommand {\T}  [1] {\bm{#1}}                        
\newcommand {\TS} [1] {\bm{#1}}                        
\newcommand {\TC} [1] {\mathbf{#1}}                    
\newcommand {\M}  [1] {\bm{#1}}                        
\newcommand {\MC} [1] {\mathbf{#1}}                    

\renewcommand {\d}            {\partial}
\newcommand   {\D}            {\mathrm{d}\mspace{1.0mu}}

\newcommand {\NM} [1] {\underaccent{\bar}{#1}}         
\newcommand {\SM} [1] {\underaccent{\tilde}{#1}}       

\providecommand {\abs} [1] {\left\vert#1\right\vert}
\providecommand {\norm}[1] {\left\Vert#1\right\Vert}

\renewcommand {\d}            {\partial}
\newcommand   {\supp}         {\operatorname{supp}}
\newcommand   {\diag}         {\operatorname{diag}}

\newcommand   {\transpose}[1] {#1^{\mathrm{t}}}
\newcommand   {\jmp}      [1] {\left\lBrack #1\right\rBrack}   
\newcommand   {\avg}      [1] {\left\lBrace #1\right\rBrace}   


\newcommand   {\ek}   {\frac{1}{2}v^2}  
\newcommand   {\ei}   {e}               
\newcommand   {\et}   {E}               
\newcommand   {\htot} {H}               

\newcommand{\N} [2][] {N_{\textsc{#2} #1}}

\title{Local high order space-time adaptive MLSDC}

\author{Erik Pfister, Matti Schulze, Jörg Stiller}
\address{%
  TU Dresden, Institute of Fluid Mechanics, 01062 Dresden, Germany}
\email{erik.pfister@tu-dresden.de}

\begin{abstract}
Building upon the semi-implicit multilevel spectral deferred correction (SI-MLSDC) method introduced in \cite{MG_Pfister2025}, this work presents a 1D space-time adaptive high-order method combining discontinuous Galerkin spectral element discretizations with multilevel spectral deferred corrections. 
The proposed approach enables genuine arbitrary-order accuracy in space and time while dynamically balancing spatial and temporal discretization errors to reduce computational cost.
A key contribution is the development of a novel temporal error estimator that in combination with a spectral error estimator in space provides a reliable basis for adaptive refinement decisions. 
The error estimator is compared with two alternative refinement criteria to assess their impact on accuracy, computational efficiency and their suitability for complex problems.
The performance of the adaptive method is demonstrated for nonlinear conservation laws ranging from Burgers' equation to the Euler equations. 
Numerical results show substantial runtime reductions while maintaining the desired accuracy. In particular, significant computational savings are achieved for Burgers' equation, and challenging benchmark problems such as the Shu-Osher shock-fluctuation benchmark.
These results demonstrate the potential of adaptive SI-MLSDC methods for efficient high-order space-time adaptive simulations of complex flow problems.
\end{abstract}


\maketitle


\section{Introduction}

The increasing use of high-order discontinuous Galerkin (DG) and spectral element discretizations for conservation laws in fluid mechanics has created a corresponding demand for time integration methods of comparable accuracy. 
This motivates the development of numerical schemes that combine arbitrarily high-order accuracy in space and time while remaining computational efficiency and stability.
To achieve high-order accuracy in time, a wide range of explicit, implicit, and semi-implicit integration methods have been proposed, including multistep and Runge-Kutta schemes \cite{TI_Gottlieb2016a,TI_Hairer1993a,TI_Hairer1996a}.
However, fully exploiting the excellent convergence properties of such methods requires an appropriate balance between spatial and temporal discretization errors.

Recent years have seen significant advances in the development of high-order semi-implicit time integration methods through the use of multi-derivative techniques \cite{TI_Frolkovic2023a,TI_Schuetz2021a,TI_Schuetz2022,TI_Zeifang2023}.
Among high-order time integration approaches, spectral deferred correction (SDC) methods provide a flexible framework for constructing arbitrarily high-order collocation schemes through iterative defect correction \cite{TI_Dutt2000a}.
Semi-implicit SDC methods naturally separate stiff and non-stiff processes and have proven particularly attractive for convection-diffusion and fluid-flow problems \cite{TI_Minion2003b,TI_Christlieb2011a,TI_Christlieb2015a,TI_Layton2005a,TI_Stiller2020a}.
Furthermore, recent developments have demonstrated stable semi-implicit SDC schemes at very high orders whith excellent stability properties \cite{TI_Stiller2024}.

A principal drawback of SDC methods is the computational cost associated with performing multiple correction sweeps.
Multilevel spectral deferred correction (MLSDC) methods address this limitation through multigrid-inspired coarse-grid corrections, thereby accelerating convergence while retaining the high-order accuracy of the underlying collocation scheme \cite{MG_Speck2014,MG_Kremling2021,MG_Hamon2019a}.
Related multilevel concepts also form the foundation of PFASST and other parallel-in-time algorithms \cite{TI_Emmett2012a,TI_Minion2010a,MG_Bolten2016a}.
More recently, robust semi-implicit MLSDC methods have been developed for conservation laws in combination with high-order DG-SEM spatial discretizations \cite{MG_Pfister2025}. This framework serves as the foundation for the developments presented in the present work.

To further improve computational efficiency, several adaptive SDC and MLSDC variants have been proposed in recent years.
Existing approaches primarily focus on adaptive time-step selection \cite{EE_Baumann2024} or adaptive mesh refinement \cite{MG_Emmett2019}.
These methods seek to preserve the high temporal accuracy of SDC while reducing the overall computational cost.
However, the corresponding studies have generally not been conducted in conjunction with genuinely high-order spatial discretizations.

Concurrently, a substantial body of literature has focused on local refinement strategies for resolving features across multiple spatial and temporal scales.
Representative examples include space-time finite element methods \cite{ST_Hansbo1994,ST_Schmich2007}, anisotropic space-time mesh adaptation \cite{ST_Jannoun2014,ST_Coupez2012}, and adaptive space-time finite element frameworks for fluid-flow problems \cite{EE_Besier2012,EE_Fidkowski2017,ST_Boisneault2023}.
Related developments also include adaptive multirate time integration methods \cite{TI_Doehring2024}.
For conservation laws, high-order ADER-DG schemes combine space-time discretizations with adaptive mesh refinement, local time stepping, and nonlinear limiting techniques \cite{SE_Dumbser2012,SE_Fambri2016,SE_Korneev2021}.
These approaches have demonstrated substantial efficiency gains for problems involving localized structures, shocks, and multiscale phenomena.

Although numerous adaptive space-time methods have been proposed, they are typically tied to specific discretization frameworks and do not provide arbitrary-order adaptivity in both space and time. To the best of the author's knowledge, such a method has not yet been investigated.
Conversely, MLSDC methods naturally provide arbitrarily high-order temporal accuracy but generally do not incorporate adaptive spatial refinement.
To the best of our knowledge, a method combining local space-time adaptivity with arbitrary-order semi-implicit MLSDC time integration has not yet been investigated.
In addition to adaptive refinement capabilities, reliable error estimation constitutes a key ingredient of any adaptive strategy.
Among the available approaches, spectral error estimators originating from \cite{EE_MavriplisDiss1989} have proven particularly attractive for high-order discretizations and have been successfully employed in adaptive DG and spectral element methods by several authors \cite{EE_Offermans2023,EE_OffermansDiss2019,EE_MassaroDiss_2024, EE_Mossier_2026}.
However, these estimators are primarily designed for spatial discretizations, and they are not directly applicable to iterative collocation methods such as SDC and MLSDC.

The present work extends the MLSDC method of \cite{MG_Pfister2025} by introducing local space-time adaptivity in order to further improve computational efficiency.
The proposed method combines dynamic space-time mesh adaptation with high-order DG-SEM discretizations and MLSDC.
In contrast to the adaptive MLSDC approach of \cite{MG_Emmett2019}, the present framework is designed for genuinely high-order discretizations in both space and time, exceeding convergence orders of four.
The resulting method supports local $h$- and $p$-refinement in space together with multilevel $p$-refinement in time, thereby enabling dynamic adaptation of both spatial and temporal approximation orders.
To guide the adaptive refinement process, a novel temporal error estimator is introduced and combined with a spatial error estimate to form a unified space-time refinement criterion.
By combining DG-SEM and MLSDC within a unified space-time multigrid method, the proposed approach retains arbitrary-order accuracy while providing local adaptivity in both space and time.
The performance of the method is assessed for scalar conservation laws, including convection-diffusion and Burgers equations, as well as for the one-dimensional Euler equations.

The remainder of this paper is organized as follows. Section~\ref{ch:problem} introduces the model problems considered throughout this work.
Section~\ref{ch:mlsdc} presents the adaptive MLSDC method, including the space-time discretization, transfer operators, adaptivity strategy, and the resulting multilevel algorithms.
In Section~\ref{ch:error_est}, a temporal error estimator is developed, a spatial error estimator is introduced, and the two are combined into a unified space-time error indicator, whose performance is subsequently evaluated.
Section~\ref{ch:tests} contains numerical experiments assessing the accuracy and efficiency of the proposed method. In addition to the error-estimator-based refinement strategy, alternative refinement criteria are investigated. The section concludes with the Shu-Osher shock-fluctuation benchmark \cite{TI_Shu1989a}.
Finally, Section~\ref{ch:conclusion} summarizes the main findings and provides concluding remarks.


\section{Problem Formulation}
\label{ch:problem}
We consider a system of conservation laws of the form
\begin{equation}
\label{eq:cons-law}
    \partial_t \mathbf{u}
    + \partial_x \mathbf{f}_c(\mathbf{u})
    =
    \partial_x \mathbf{f}_d(\mathbf{u})
    + \mathbf{f}_s(\mathbf{u},x,t),
\end{equation}
where $\mathbf{u}(x,t)\in\mathbb{R}^d$ denotes the vector of conserved variables, $x\in\Omega\subset\mathbb{R}$ is the spatial coordinate, and $t\in\mathbb{R}^{+}$ is time. 
The functions $\mathbf{f}_c$, $\mathbf{f}_d$, and $\mathbf{f}_s$ represent the convective fluxes, diffusive fluxes, and source terms, respectively. Introducing the combined right-hand side operator
\begin{equation}
    \mathbf{f}
    =
    -\partial_x \mathbf{f}_c
    + \partial_x \mathbf{f}_d
    + \mathbf{f}_s,
\end{equation}
the governing equations can be written in evolutionary form as
\begin{equation}
    \partial_t \mathbf{u}
    =
    \mathbf{f}(\mathbf{u},x,t).
\end{equation}
Let $\mathbf{A}_c = \mathbf{f}_c'(\mathbf{u})$ denote the Jacobian of the convective flux. Assuming a diffusive flux of the form
\begin{equation}
    \mathbf{f}_d
    =
    \mathbf{A}_d \, \partial_x \mathbf{u},
\end{equation}
where $\mathbf{A}_d(\mathbf{u},\partial_x\mathbf{u})$ is the diffusion matrix, the system can be expressed in quasilinear form as
\begin{equation}
    \partial_t \mathbf{u}
    =
    -\mathbf{A}_c \, \partial_x \mathbf{u}
    + \partial_x \left( \mathbf{A}_d \, \partial_x \mathbf{u} \right)
    + \mathbf{f}_s.
\end{equation}

\subsection{Scalar conservation laws}
As representative scalar model problems, we consider the case where the solution vector $\mathbf{u}$ reduces to a single scalar variable $u$. 
In particular, we study the convection-diffusion equation
\begin{align}
\label{eq:conv-diff}
\partial_t u
=
-\partial_x (v u)
+
\partial_x \left( \nu \partial_x u \right),
\end{align}
where $v$ denotes the advection velocity and $\nu \ge 0$ the diffusion coefficient. 
Furthermore, we consider the viscous Burgers equation with constant diffusivity,
\begin{align}
\label{eq:burgers}
\partial_t u
=
-\partial_x \left( \tfrac{1}{2} u^2 \right)
+
\partial_x \left( \nu \partial_x u \right)
+
f_{\mathrm{s}}(x,t),
\end{align}
where $f_{\mathrm{s}}$ is a prescribed source term.

\subsection{Euler and  Navier-Stokes equations}
\label{sec:conservation-laws:cns}
As a more complex test case, we consider the one-dimensional Euler and compressible Navier-Stokes equations. The vector of conserved variables is given by
\begin{equation}
    \mathbf{u}
    =
    \begin{bmatrix}
        \rho \\
        \rho v \\
        \rho e_t
    \end{bmatrix},
\end{equation}
where $\rho$ denotes the density, $v$ the velocity, and $e_t$ the specific total energy. 
The corresponding convective and diffusive fluxes are
\begin{equation}
    \mathbf{f}_c
    =
    \begin{bmatrix}
        \rho v \\
        \rho v^2 + p \\
        \rho v e_t + p v
    \end{bmatrix},
    \qquad
    \mathbf{f}_d
    =
    \begin{bmatrix}
        0 \\
        \frac{4}{3}\eta\,\partial_x v \\
        \frac{4}{3}\eta\,v\,\partial_x v + \lambda\,\partial_x T
    \end{bmatrix},
\end{equation}
where $p$ is the pressure, $T$ the temperature, $\eta$ the dynamic viscosity, and $\lambda$ the thermal conductivity.
The convective Jacobian $\mathbf{A}_c=\mathbf{f}_c'(\mathbf{u})$ is given in Eq.~(8) of \cite{MG_Pfister2025}. Likewise, the diffusion matrix $\mathbf{A}_d$ is provided in Eq.~(9) of \cite{MG_Pfister2025}.


\newcommand{\Restrict}    {\NM{\mathcal R}}
\newcommand{\Project}     {\NM{\mathcal P}}
\newcommand{\Interpolate} {\NM{\mathcal I}}

\section{Adaptive MLSDC}
\label{ch:mlsdc}
\subsection{Space-time discretization and collocation method}
\label{sec:discretization}
In preparation for the adaptive multilevel method, we first introduce the notation employed throughout this paper. 
The combination of spatial and temporal discretizations, SDC and multiple levels leads to a relatively rich indexing structure. To facilitate readability, the most important symbols and conventions are summarized in Table~\ref{tab:notation}. 
Boldface quantities denote arrays containing all solution components, as shown in Section \ref{ch:problem}. 
A sequence of levels is considered and indexed by
$l=1,\ldots,L$.
Level $l=1$ denotes the coarsest level, whereas $l=L$ denotes the finest.
\begin{table}[H]
\centering
\caption{Notation used for continuous and discrete solution representations and space-time discretizations.}
\label{tab:notation}
\small
\begin{tabular}{p{0.22\textwidth} p{0.55\textwidth} p{0.15\textwidth}}
\toprule
Notation & Description & First Ref. \\
\midrule
$\mathbf{u}(x,t)$
    & Continuous solution
    & Eq.~\eqref{eq:cons-law} \\
$\mathbf{u}_{h,l}(x,t)$
    & Approximate continuous solution on level $l$
    & Eq.~\eqref{eq:approx_solution} \\
$\mathbf{u}_{m,l}(x)$
    & Approximate solution at temporal collocation node $m$ on level $l$
    & Eq.~\eqref{eq:solution_node_m} \\
$\mathbf{u}_{m,l}^k(x)$
    & Approximate solution at temporal collocation node $m$ on level $l$ at sweep $k$
    & Eq.~\eqref{eq:SL_corrector} \\
$\NM{\mathbf{u}}_{m,l}$
    & Spatial nodal coefficients at temporal collocation node $m$ on level $l$
    & Eq.~\eqref{eq:spatial_coeffs} \\
$\NM{\mathbf{U}}_l$
    & Spatio-temporal nodal coefficient vector on level $l$
    & Eq.~\eqref{eq:spatio-temporal_coeffs} \\
$\NM{\mathbf{F}}(\NM{\mathbf{U}})$
    & Collocation operator
    & Eq.~\eqref{eq:cm_incremental_matrix} \\
$\NM{\mathbf{f}}(\NM{\mathbf{u}})$
    & Right-hand side of the collocation method
    & Eq.~\eqref{eq:cm_incremental} \\
$\NM{M}_l^x$
    & Spatial mass matrix on level $l$
    & Eq.~\eqref{eq:cm_incremental} \\
$\mathbf{u}^{e,n}_{i,m,l}$
    & Nodal DG coefficient in element $e$, time step $n$, spatial node $i$,
      collocation node $m$, on level $l$
    & Eq.~\eqref{eq:space-time_ansatz} \\
$T_l$
    & Temporal domain on level $l$
    & Eq.~\eqref{eq:domain} \\
$\Omega_l$
    & Spatial domain on level $l$
    & Eq.~\eqref{eq:domain} \\
$T_l^n$
    & Time step on level $l$
    & Eq.~\eqref{eq:space-time_elements} \\
$\Omega_l^e$
    & Spectral element $e$ on level $l$
    & Eq.~\eqref{eq:space-time_elements} \\
$Q_l^{e,n}$
    & Space-time element on level $l$
    & Eq.~\eqref{eq:space-time_elements} \\
\bottomrule
\end{tabular}
\end{table}
We now consider a hierarchy of space-time discretizations
\begin{align}
\label{eq:domain}
Q_l = \Omega_l \times T_l ,
\end{align}
where $\Omega_l$ denotes the spatial discretization on level $l$. 
On each level $l$, the temporal domain $T_l$ is partitioned into time intervals
\begin{align}
T_l^n=[t^{n-1},t^n]
\end{align}
of uniform length $\Delta t_l = \Delta t$. Within each time interval, a set of temporal collocation nodes is introduced,
\begin{align}
t^{n-1}
=:\ t^n_{0,l} < t_{1,l} < \dots < t_{m,l} < \dots < t^n_{M_l,l}
=
t^{n}.
\end{align}
The number of collocation nodes $M_l$ may vary between levels.
Different choices of collocation nodes lead to different properties of the resulting time-integration scheme. Unless stated otherwise, Radau-Right nodes are employed throughout this work. 
The semidiscrete solution at the collocation nodes is represented by
\begin{align}
    \label{eq:solution_node_m}
    \mathbf{u}_{m,l} \approx \mathbf{u}(x,t_{m,l}).
\end{align}
\begin{figure}[H]
    \centering
    \includegraphics[width=0.55\linewidth]{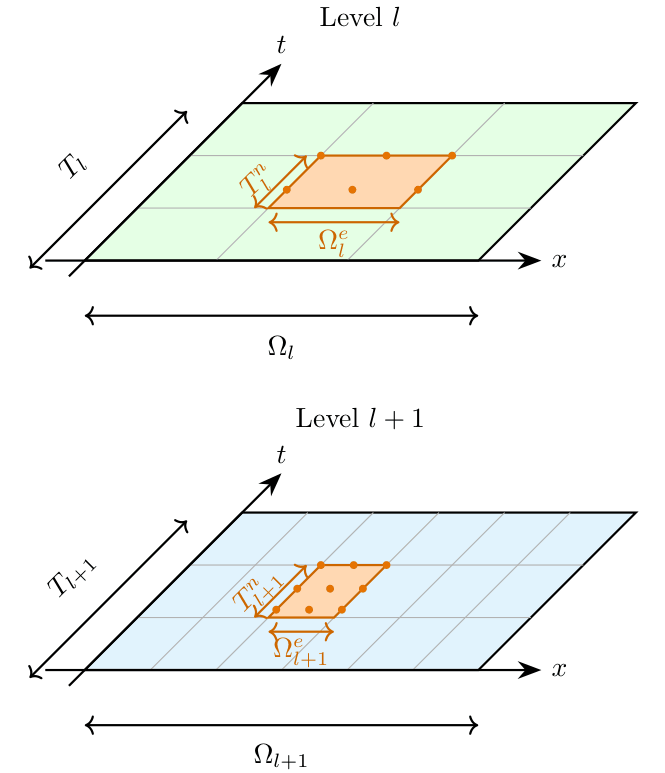}
    \caption{Sketch of the multilevel space-time discretization for two levels. Level $l+1$ exemplary employs additional spatial elements ($h$-refinement) and more collocation nodes in time ($p$-refinement) than level $l$. Orange points indicate the space-time coefficients within the space-time elements $Q_l^{e,n}$.}
    \label{fig:space-time_discr}
\end{figure}
The spatial discretization is based on the DG-SEM, in which the physical domain $\Omega_l$ is partitioned into a set of non-overlapping spectral elements $\Omega_l^e$. The computational domain on level $l$ is therefore given by
\begin{align}
\label{eq:elements}
\Omega_l = \bigcup_e \Omega_l^e .
\end{align}
The combined space-time elements are denoted by
\begin{align}
\label{eq:space-time_elements}
Q_l^{e,n} = \Omega_l^e \times T_l^n .
\end{align}
An illustration of the space-time discretized domain is shown in Figure~\ref{fig:space-time_discr}.
On each level $l$, the approximate solution $\mathbf{u}_{h,l}(x,t)$ is sought in a finite-dimensional space of piecewise polynomial functions,
\begin{align}
\label{eq:approx_solution}
\mathbb{U}_{h,l}
=
\left\{
    \mathbf{u}_{h,l} \in \left[\mathbb{L}^2(Q_l)\right]^d :
    \mathbf{u}_{h,l}\big|_{Q_l^{e,n}}
    \in
    \left[
        \mathbb{P}_{P_l}(\Omega_l^e)
        \otimes
        \mathbb{P}_{M_l-1}(T_l^n)
    \right]^d,
    \quad
    \forall\, Q_l^{e,n} \subset Q_l
\right\},
\end{align}
The polynomial space-time ansatz on level $l$ for a space-time element $Q_l^{e,n}$ is given by
\begin{align}
\label{eq:space-time_ansatz}
\mathbf{u}^{e,n}_{h,l}(\xi,\tau)
&=
\sum_{i=0}^{P_l}
\sum_{m=1}^{M_l}
\mathbf{u}^{e,n}_{i,m,l}
\, \ell_{i,l}^x(\xi)\,\ell_{m,l}^t(\tau).
\end{align}
The functions $\ell^x_{i,l}(\xi)$ and $\ell^t_{m,l}(\tau)$ denote the spatial and temporal nodal Lagrange basis functions, respectively. Here, $\xi = \xi^e(x) \in[-1,1]$ and $\tau = \tau^n(t) \in[0,1]$ are the standard spatial and temporal  reference coordinates associated with the space-time element $Q_l^{e,n}$.
The nodal spatio-temporal coefficients $\mathbf{u}^{e,n}_{i,m,l}$ define the space-time approximation on the element $Q_l^{e,n}$. 
The spatial coefficients associated with collocation node $m$ are grouped into
\begin{align}
\label{eq:spatial_coeffs}
\NM{\mathbf{u}}_{m,l}
=
\bigl[
\{
\mathbf{u}^{e}_{i,m,l}
\}_{i=0}^{P_l}
\bigr].
\end{align}
The collection of all collocation-node coefficient vectors on level $l$ is denoted by
\begin{align}
\label{eq:spatio-temporal_coeffs}
\NM{\mathbf{U}}_l
=
\bigl[
\{
\NM{\mathbf{u}}_{m,l}
\}_{m=1}^{M_l}
\bigr].
\end{align}
Having introduced the underlying space-time discretization, we now derive the collocation problem underlying the SDC method.
Starting from the Picard integral formulation of a generic initial value problem, a collocation discretization is obtained by applying a quadrature rule over successive subintervals, yielding
\begin{equation}
\label{eq:cm_incremental}
  \NM{M}_l^x
  \NM{\mathbf u}_{m,l}
  =
  \NM{M}_l^x
  \NM{\mathbf u}_{m-1,l}
  +
  \Delta t
  \sum_{j=1}^{M_l}
  w^{\mathsc{nn}}_{m,j,l}
  \NM{\mathbf f}(\NM{\mathbf u}_{j,l}),
  \quad
  m = 1,\dots,M_l .
\end{equation}
Here, $w^{\mathsc{nn}}_{m,j,l}$ denote the incremental integration weights, defined by
\begin{align}
\label{eq:nn_weights}
w^{\mathsc{nn}}_{m,j,l}
=
\frac{1}{\Delta t}
\int_{\tau_{m-1,l}}^{\tau_{m,l}}
\ell_{j,l}^t(t)\,\D \tau .
\end{align}
Equation~\eqref{eq:cm_incremental} defines the \emph{node-to-node} (or incremental) collocation method (CM). 
Please note: the spatial mass matrix $\NM{M}_l^x$ should not be confused with the number of temporal collocation nodes $M_l$. 
In this work, we exclusively consider the incremental formulation, the alternative non-incremental form and the reasons for not employing it are discussed in \cite{MG_Pfister2025}. 
The collocation method may be interpreted as a system of equations for the solution at all collocation nodes:
\begin{equation}
\label{eq:cm_incremental_matrix}
\NM{\mathbf F}_l(\NM{\mathbf U}_l)
=
\mathbf{ \NM 0} .
\end{equation}
Further details on the temporal and spatial discretizations are available in \cite{MG_Pfister2025,TI_Stiller2024}. The fully discrete collocation formulation is not repeated here and can be found in Eq.~(56) of \cite{MG_Pfister2025}.

\subsection{Local adaptivity strategy}
\label{sec:mlsdc_markers}
\begin{figure}[H]
    \centering
    \includegraphics[width=0.85\linewidth]{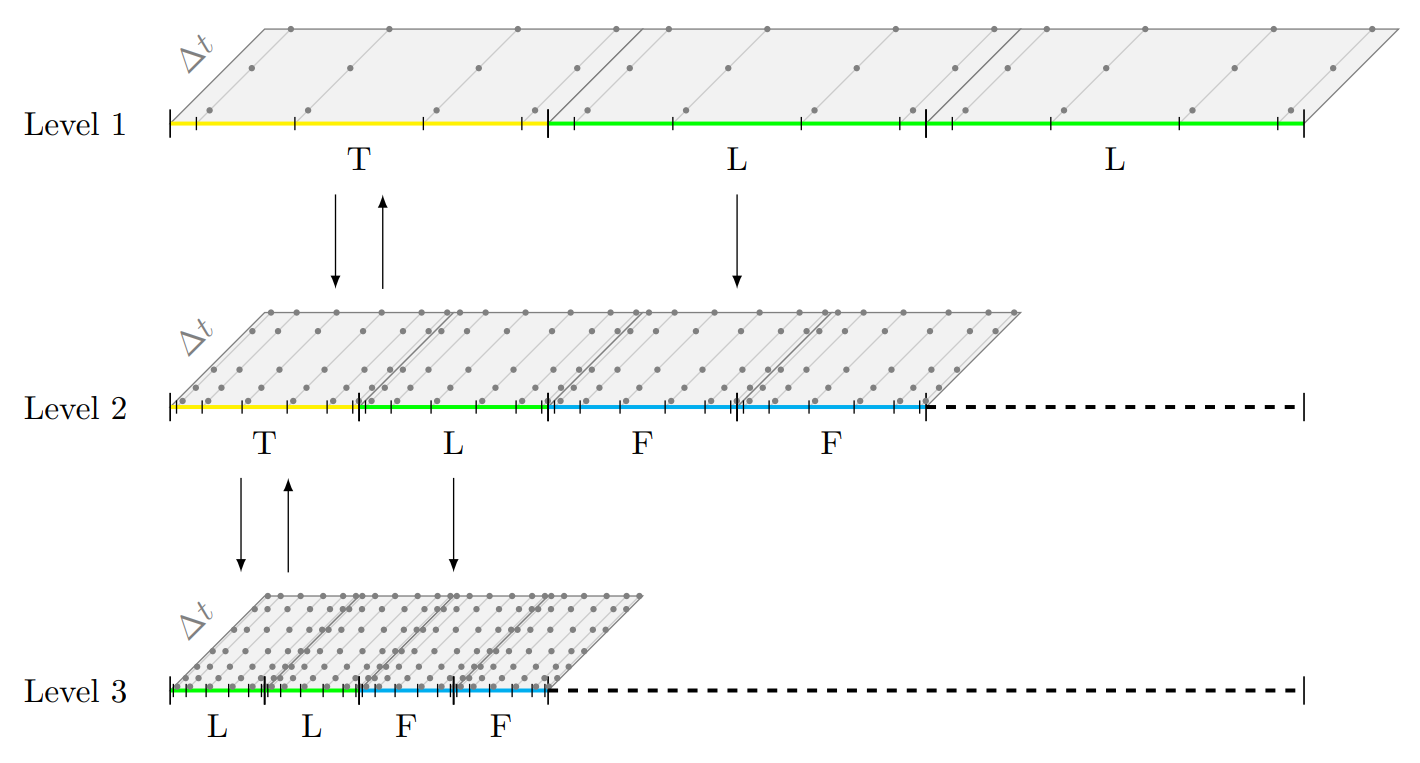}
    \caption{Sketch of the adaptivity routine, showing $h/p$-refinement in space and $p$-refinement in time for the first refinement and $h$-refinement in space and $p$-refinement in time for the second refinement. The example contains three elements on the coarsest level, non-periodic boundary conditions and increasing numbers of collocation nodes across levels, e.g. $M_l=\{3,5,7\}$. $\mathsc{t}$ - twig, $\mathsc{l}$ - leaf, $\mathsc{f}$ - frozen.}
    \label{fig:adaptivity}
\end{figure}

The preceding discussion considered a globally refined multilevel discretization and its associated CM formulation.
However, the present work employs an adaptive MLSDC method, requiring a distinction between different element classes in the multilevel hierarchy, as sketched in Figure~\ref{fig:adaptivity}. 
Throughout this work, the superscripts $\mathsc{a}$, $\mathsc{f}$, $\mathsc{t}$, and $\mathsc{l}$ denote active, frozen, twig, and leaf elements, respectively. For example, the active region on level $l$ is defined as
\begin{align}
\Omega_l^\mathsc{a}
=
\Omega_l^\mathsc{t}
\cup
\Omega_l^\mathsc{l}
=
\Omega_l^\mathsc{tl}
=
\Omega_l \setminus \Omega_l^\mathsc{f}.
\end{align}
Leaf elements are not refined further and therefore represent the terminal elements of the adaptive hierarchy. At inner-subdomain boundaries on the next finer level, they provide interpolated values for frozen elements. Twig elements, in contrast, may be refined and serve as parents of leaf elements on the finer level.
The adaptive MLSDC method extends a globally refined multilevel formulation by restricting computations to regions where additional resolution is required. The final goal is enabling local space-time adaptivity. 
Unlike the globally refined approach of \cite{MG_Pfister2025}, computations are performed only on elements that actively participate in the multilevel hierarchy.
The framework supports both $h$- and $p$-refinement in space, whereas refinement in time is restricted to increasing the number of collocation nodes across levels, hence $p$-refinement.
Each element is assigned an activity marker and a refinement marker defining its element class. 
Active elements comprise both twig and leaf elements. The refinement and activity markers are controlled by the error estimator described in Section~\ref{ch:error_est}. On the coarsest level, all elements are active.
Calculations are performed exclusively on active elements. Values stored on frozen elements are interpolated from coarser parent levels (leafs) and provide boundary data for adjacent active regions, see Figure~\ref{fig:adaptivity}. 
For the initial values at the left-most temporal node ($m=0$), interpolation from a coarser level is performed only when a subdomain is newly refined. 
The resulting hierarchy consists of locally refined space-time regions embedded within coarser surrounding domains while retaining the multilevel structure required by the MLSDC algorithm.
The active region is enveloped by a layer of frozen elements that do not participate in the multilevel iteration. Nevertheless, they retain the standard DG element structure and are included in the evaluation of numerical fluxes and element coupling terms. Consequently, interfaces between active and frozen elements are treated identically to ordinary inter-element interfaces. This avoids the introduction of artificial Dirichlet boundary conditions on inner-subdomain boundaries and preserves the original DG discretization locally.

\subsection{Transfer operators}
\label{sec:transferoperators}
The multilevel formulation requires transfer operators to exchange information between adjacent levels. 
These operators are represented by algebraic matrices acting on the corresponding coefficient vectors. To keep the notation compact, each transfer operator is identified only by the index of its target level, while the source level is implied by the transfer direction.
For notational purposes, the coefficient vector on level $l$ is decomposed according to the element classes,
\begin{equation}
    \NM{\mathbf U}_l
    =
    \begin{bmatrix}
        \NM{\mathbf U}_l^{\mathsc f} \\
        \NM{\mathbf U}_l^{\mathsc t} \\
        \NM{\mathbf U}_l^{\mathsc l}
    \end{bmatrix},
\end{equation}
where the subvectors correspond to frozen, twig, and leaf elements, respectively. 
The active coefficient vector is given by
\begin{equation}
    \NM{\mathbf U}_l^{\mathsc a}
    =
    \begin{bmatrix}
        \NM{\mathbf U}_l^{\mathsc t} \\
        \NM{\mathbf U}_l^{\mathsc l}
    \end{bmatrix},
\end{equation}
and collects all degrees of freedom that participate in the adaptive multilevel iteration.
\begin{remark}
The transfer operators are formally defined on the complete coefficient vectors. In the adaptive MLSDC algorithm, however, their action is restricted by the adaptive multilevel hierarchy. \emph{Interpolation} transfers information from coarse elements to finer elements within the active region and additionally facilitates the construction of frozen elements. \emph{Restriction and projection} collect information from active finer child elements and transfer it to the corresponding twig elements on the next coarser level.
\end{remark}
Coarse-to-fine transfer of solution values and corrections is performed by the interpolation operator,
\begin{equation}
    \NM{\mathbf U}_{l}
    =
    \mathcal I_l
    \NM{\mathbf U}_{l-1},
\end{equation}
which employs embedded interpolation such that coarse-grid information is represented on the finer level up to round-off accuracy.
Fine-to-coarse transfer of residuals is carried out by the restriction operator,
\begin{equation}
    \NM{\mathbf R}_{l}
    =
    \mathcal R_l
    \NM{\mathbf R}_{l+1},
\end{equation}
which aggregates residual information from active child elements and transfers it to the corresponding twig elements on the next coarser level.
In the present work, the restriction operator is chosen as the transpose of the interpolation operator,
\begin{equation}
    \mathcal R_l
    =
    \bigl(
        \mathcal I_{l+1}
    \bigr)^t,
\end{equation}
although alternative definitions, such as an $L^2$-projection-based restriction, are also possible.
The full approximation scheme (FAS) additionally requires a fine-to-coarse projection operator for solution values,
\begin{equation}
    \NM{\mathbf U}_{l}
    =
    \mathcal P_l
    \NM{\mathbf U}_{l+1}.
\end{equation}
Unlike residual restriction, the definition of a suitable solution projection is not straight forward.
Within the present high-order DG framework, two projection variants are considered: embedded interpolation and the $L^2$-projection. A special situation arises for $2{:}1$ $h$-coarsening, where two neighboring fine elements are merged into a single coarse element and the projection operator must consistently combine information from both finer child elements.
Unless stated otherwise, the embedded interpolation as projection is used for tests.
For embedded interpolation, each fine element contributes only to the corresponding portion of the coarse element. 
Double contributions due to jumps are averaged.
In contrast, the $L^2$-projection incorporates information from both child elements into the complete coarse-element representation.
For both restriction and projection, only active child elements contribute to the transfer, whereas the resulting quantities are stored on the corresponding twig elements of the next coarser level. 
In locally adaptive settings, transfer operators therefore act only on active regions. Values on frozen elements are obtained by interpolation. 
Spatial and temporal transfer operators are constructed analogously, since both rely on polynomial approximation spaces. The corresponding transfers are applied successively and therefore act together as coupled space-time operators.
Optional modifications may be introduced in the vicinity of discontinuities. In the present work, these include linear blending as the default strategy and coefficient averaging at element interfaces.

\subsection{Multilevel collocation}
\label{sec:ml_collocation}
Brandt's FAS multigrid work provides the theoretical foundation for the present formulation, a comprehensive derivation can be found in \cite{MG_Brandt2011}. 
Further details on the underlying multilevel collocation or MLSDC method and its  formulation are given in \cite{MG_Pfister2025}.
For globally refined meshes, the starting point for the multilevel collocation problem is the incremental CM \eqref{eq:cm_incremental_matrix} and may be written as
\begin{equation}
\label{eq:cm_ml}
    \NM{\mathbf F}_l(\NM{\mathbf U}_l)
    =
    \NM{\mathbf G}_l,
    \qquad
    l = 1,\dots,L,
\end{equation}
where the additional term FAS right-hand side $\NM{\mathbf G}_l$ incorporates information from finer grids through the FAS correction to ensure that all levels approach the finest one.
The central idea of the locally adaptive MLSDC method is to construct and update the FAS right-hand side only in regions where refinement is required. Consequently, the multilevel correction process is restricted to active portions of the hierarchy, while inactive regions are excluded from the full multilevel iteration.
Considering an approximate solution $\NM{\mathbf{\tilde U}}_l$, the locally adaptive FAS-DG formulation reads
    \begin{subequations}
    \label{eq:helmholtz:fas}
    \begin{alignat}{2}
    \label{eq:helmholtz:fas:uf}
    \NM{\mathbf{\tilde U}}_l^{\mathsc f}
    &=
    \mathcal I_l^{\mathsc{f}}
    \NM{\mathbf{\tilde U}}_{l-1}
    \qquad
    &&\text{in } \Omega_l^{\mathsc f}
    \\[5pt]
    \label{eq:helmholtz:fas:system}
    \NM{\mathbf F}_l^{\mathsc a}
    \!\left(
    \NM{\mathbf{\tilde U}}_l
    \right)
    &=
    \NM{\mathbf{ G}}_l^{\mathsc a}
    \qquad
    &&\text{in } \Omega_l^{\mathsc a}
    \\[5pt]
    \label{eq:helmholtz:fas:rhs}
    \NM{\mathbf{ G}}_l^{\mathsc a}
    &=
    \begin{bmatrix}
    \NM{\mathbf{ G}}_l^{\mathsc t}
    \\[3pt]
    \mathbf{ \NM 0}
    \end{bmatrix}
    \qquad
    &&\text{in } \Omega_l^{\mathsc a}
    \\[5pt]
    \label{eq:helmholtz:fas:rhs_twig}
    \NM{\mathbf{ G}}_l^{\mathsc t}
    &=
    \NM{\mathbf F}_l^{\mathsc t}
    \!\left(
    \mathcal P_l^{\mathsc{t}}
    \NM{\mathbf{\tilde U}}_{l+1}
    \right)
    +
    \mathcal R_l^{\mathsc{t}}
    \NM{\mathbf{ R}}_{l+1},
    \qquad l < L
    \\[5pt]
    \label{eq:helmholtz:fas:residual}
    \NM{\mathbf{ R}}_l^{\mathsc a}
    &=
    \NM{\mathbf{ G}}_l^{\mathsc a}
    -
    \NM{\mathbf F}_l^{\mathsc a}
    \!\left(
    \NM{\mathbf{\tilde U}}_l
    \right)
    \qquad
    &&\text{in } \Omega_l^{\mathsc a}.
    \end{alignat}
    \end{subequations}
where
$\NM{\mathbf{ R}}_l^{\mathsc a}$
denotes the active residual.
Equation~\eqref{eq:helmholtz:fas:system} defines the active coefficients
$\NM{\mathbf{\tilde U}}_l^{\mathsc a}$ through the active subset of collocation equations and the active right-hand side
$\NM{\mathbf{ G}}_l^{\mathsc a}$.
The FAS right-hand side is required only on levels $l<L$. On the finest level $L$, no finer approximation exists and consequently
\begin{equation*}
\NM{\mathbf{ G}}_L
=
\NM{\mathbf 0}.
\end{equation*}
The composition of the modified right-hand side depends on the element class. For leaf elements it vanishes, such that the original collocation equations are recovered. For twig elements it consists of the collocation operator evaluated using projected child information together with the restricted child residual.
Equation~\eqref{eq:helmholtz:fas:residual} defines the residual as the difference between the modified right-hand side and the active collocation equations.
Upon convergence, the residual vanishes and the adaptive FAS formulation reduces to
\begin{subequations}
\label{eq:helmholtz:fas:converged}
\begin{align}
\label{eq:helmholtz:fas:converged:frozen}
\NM{\mathbf U}_l^{\mathsc f}
&=
\mathcal I_l^{\mathsc{f}}
\NM{\mathbf U}_{l-1}
\\
\label{eq:helmholtz:fas:converged:twigs}
\NM{\mathbf U}_l^{\mathsc t}
&=
\mathcal P_l^{\mathsc{t}}
\NM{\mathbf U}_{l+1}
\\
\label{eq:helmholtz:fas:converged:leaves}
\NM{\mathbf F}_l^{\mathsc l}
\!\left(
\NM{\mathbf U}_l
\right)
&=
\NM{\mathbf 0}.
\end{align}
\end{subequations}
It is important to emphasize that both the adaptive formulation \eqref{eq:helmholtz:fas} and its converged form \eqref{eq:helmholtz:fas:converged} remain fully consistent with the underlying DG discretization. Consequently, the conservation properties of the original scheme are preserved on every level of the adaptive hierarchy.
\begin{remark}
The globally refined MLSDC method solves to the collocation problem on the finest level (see, e.g. \cite{MG_Pfister2025}) and inherits its properties. In the present setting with Radau Right (RR) nodes, this collocation problem corresponds to the Radau IIA method. 
\end{remark}
\begin{remark}
Conservation properties are satisfied independently on each level. However, the coupling between levels through the transfer operators and the FAS right-hand side contribution $\NM{\M G}_l$, acting as an additional source term, may introduce non-conservative information transfer between consecutive levels.
\end{remark}

\subsection{Multilevel SDC algorithm}
The SDC method serves as an iterative solver for the collocation problem by successively applying corrections based on low-order time-integration methods.
Detailed derivations of the SDC framework can be found in \cite{TI_Dutt2000a,TI_Minion2003b,TI_Ong2020a,TI_Stiller2020a}. 
Since SDC is an iterative method, the approximation at iteration $k$ is denoted by $\mathbf{u}^{k}_{m,l}(x)$, where a single iteration is commonly referred to as a \emph{sweep}.
One- or two-staged predictors and correctors are used for the SDC scheme in the present manuscript.
The predictor employs a low-order time-stepping scheme to advance through the collocation nodes and generate an initial approximation $\mathbf{u}^{0}_{m,l}$. 
In this work, the predictor and corrector are based on IMEX Euler as well as the one- and two-stage semi-implicit methods SI(1) and SI(2) introduced in \cite{TI_Stiller2024}. Further implementation details are provided in Section~3.3 of \cite{MG_Pfister2025}.
In the multilevel context of MLSDC there is an additional term for the corrector equation - the FAS-RHS $\NM{\mathbf G}_l$.
The MLSDC correction equation for the $k$-th approximation is given by
\begin{equation}
\label{eq:SL_corrector}
  \NM{M}_l^x
  \NM{\mathbf u}^{k,\mathsc{A}}_{m,l}
  =
  \NM{M}_l^x
  \NM{\mathbf u}^{k,\mathsc{A}}_{m-1,l}
  +
  \Delta t_l
  \sum_{j=1}^{M_l}
  w^{\mathsc{nn}}_{m,j,l}
  \NM{\mathbf f}(\NM{\mathbf u}^{k-1,\mathsc{A}}_{j,l})
  +
  \Delta \NM{\mathbf H}^{k,\mathsc{A}}_{m,l}
  + \NM{\mathbf g}_{m,l}^\mathsc{A},
  \quad
  m = 1,\dots,M_l .
\end{equation}
Here, the term $\Delta \NM{\mathbf H}^{k}_{m,l}$ characterizes the underlying SDC method, such as IMEX Euler, SI(1), or SI(2) and results from a semi-implicit time integration over $T_m^n$. Its precise definition, together with additional implementation details, is provided in Section~3.3 of \cite{MG_Pfister2025}. 
As the method converges, the correction term vanishes,
\[
\Delta \NM{\mathbf H}^{k,\mathsc{A}}_{m,l} \rightarrow \mathbf{\NM 0},
\]
and Eq.~\eqref{eq:SL_corrector} recovers the collocation formulation (Section~\ref{sec:ml_collocation}). 
The use of Radau-Right collocation nodes yields the Radau IIA method of order $2M_l-1$, see \cite[Ch.~6]{TI_Deuflhard2002a}.
The residual of the incremental MLSDC method for the $k$-th approximation reads
\begin{equation}
\label{eq:mlsdc:incremental:residual}
  \NM{\mathbf r}^{k,\mathsc{A}}_{m,l}
  =
  \NM{\mathbf g}^{\mathsc{A}}_{m,l}
  -
  \NM{M}_l^x
  \NM{\mathbf u}^{k,\mathsc{A}}_{m,l}
  +
  \NM{M}_l^x
  \NM{\mathbf u}^{k,\mathsc{A}}_{m-1,l}
  +
  \Delta t_l
  \sum_{j=1}^{M_l}
  w^{\mathsc{nn}}_{m,j,l}
  \NM{\mathbf f}(\NM{\mathbf u}^{k-1,\mathsc{A}}_{j,l}),
  \quad
  m = 1,\dots,M_l .
\end{equation}
Please note: for reasons of readability, this convention applies to the following Algorithm \ref{alg:corrector_mlsdc}: wherever no superscript index appears next to the required quantities, active elements are meant, also the index $l$ is dropped for a more concise style.
\begin{algorithm}[H]
\caption{Locally adaptive incremental MLSDC correction sweeps (for active elements).}
\label{alg:corrector_mlsdc}
\begin{algorithmic}[1]
\Procedure{MLSDC\_Corrector}{}${(\NM{\M U},\NM{\M G}, N)}$
  \State
    ${\NM{\M U}^{0} \gets \NM{\M U}}$
  \For{$k = 1,N$}
    \State
      ${\NM{\M U}_0^k \gets \NM{\M U}_0}$
    \For{$m = 1,M$}
            \State
                ${\NM{\M u}_{m}^{k} \gets
                \textsc{Solve}
                    \big( \NM{M}^x \NM{\M u}^{k}_{m}
                        = \NM{M}^x \NM{\M u}^{k}_{m-1}
                        + \Delta t \sum_{j=1}^{M} w_{m,j} \M f(\NM{\M u}^{k-1}_{j},t_{j})
                        + \Delta \NM{\M H}^{k}_{m} + \NM{\M g}_{m}
                    \big)
                }$
    \EndFor
  \EndFor
  \State
    ${\NM{\M U} \gets \NM{\M U}^{N}}$
\EndProcedure
\end{algorithmic}
\end{algorithm}
The adaptive MLSDC method employs the same V-cycle structure as described in \cite{MG_Pfister2025}, consisting of pre-smoothing, restriction, coarse-level solution, correction prolongation, and, optionally, post-smoothing.
To account for local mesh refinement, the V-cycle is modified such that FAS right-hand sides are constructed only on active regions and corrections are applied exclusively to active elements (see Algorithm~\ref{alg:mlsdc:v-cycle}) as determined by the activity and refinement markers.
Since the adaptive hierarchy changes over time, suitable starting strategies are required to initialize the solution on newly activated fine-grid elements prior to each time step.
The Cascade, FMG1, and FMG2 initialization strategies are employed for this purpose, a detailed description and discussion of these approaches is given in \cite{MG_Pfister2025}.
For the FMG-based initialization strategies, the auxiliary V-cycles are likewise restricted to the active portions of the hierarchy, thereby maintaining consistency with the adaptive multilevel formulation.
\begin{algorithm}[H]
\caption{Locally adaptive MLSDC V-cycle.}
\label{alg:mlsdc:v-cycle}
\begin{algorithmic}[1]
\Procedure{MLSDC\_V\_Cycle}{}{$(\NM{\M U}_{1:L}, N_{c}, L)$}
  \State
    ${\NM{\M G} \; \gets \; \NM{\M 0}}$
  \State 
    $\NM{\M R} \; \gets \; \NM{\M 0}$
  \For{$l = L,2,-1$}
    \State {${\NM{\M U}_{l}^{\mathsc{A}}}$}
      $ {\gets \; {\textsc{MLSDC\_Corrector}(\NM{\M U}_l, \NM{\M G}_l, N_s = 1)}}$
      \textcolor{Grey}{\Comment{{Pre-smoothing}}}
    \State
      \makebox[2.2em][l]{${\NM{\M R}_{l}^{\mathsc{A}}}$}
        ${\gets \; \NM{\M G}_l^{\mathsc{A}} -\NM{\M F}_{l}^{\mathsc{A}}(\NM{\M U}_{l}) }$
        \textcolor{Grey}{\Comment{{Residual evaluation}}}
    \State
      \makebox[2.2em][l]{${\NM{\M V}_{l-1}^{\mathsc{T}}}$}
        ${ \gets \; \mathcal{P}_{l-1}^{\mathsc{A}}\NM{\M U}_{l}}$
        \textcolor{Grey}{\Comment{{Solution restriction}}}
        \State
          \makebox[2.2em][l]{${\NM{\M G}_{l-1}^{\mathsc{T}}}$}
            ${ \gets \; \NM{\M F}_{l-1}^{\mathsc{T}}(\NM{\M V}_{l-1})
                      + \mathcal{R}_{l-1}^{\mathsc{T}}\NM{\M R}_{l}}
            $
          \textcolor{Grey}{\Comment{{RHS composition}}}
  \EndFor
  \State {${\NM{\M U}_{1}^{\mathsc{A}}}$}
      $ {\gets \; {\textsc{MLSDC\_Corrector}(\NM{\M U}_1, \NM{\M G}_1, N_{c})}}$
    \textcolor{Grey}{\Comment{{Coarse solution}}}
  \For{$l = 2,L$}
      \State
        \makebox[1.5em][l]{${\NM{\M U}_{l}^{\mathsc{A}}}$}
          ${ \gets \; \NM{\M U}_{l}^{\mathsc{A}}
                    + \Interpolate_{l}^{\mathsc{A}}(\NM{\M U}_{l-1} - \NM{\M V}_{l-1})}$
          \textcolor{Grey}{\Comment{{Correction prolongation}}}
      \State
        \makebox[1.5em][l]{${\NM{\M U}^{\mathsc{f}}_{l}}$}
          ${ \gets \; \Interpolate^{\mathsc{f}}_{l}\NM{\M U}_{l-1}}$
          \textcolor{Grey}{\Comment{{Interpolate boundary conditions}}}
    \If{$l < L$}
      \State {${\NM{\M U}_{l}^{\mathsc{A}}}$}
      $ {\gets \; {\textsc{MLSDC\_Corrector}(\NM{\M U}_l, \NM{\M G}_l, N_s = 1)}}$
      \textcolor{Grey}{\Comment{{Post-smoothing}}}
    \EndIf
  \EndFor
\EndProcedure
\end{algorithmic}
\end{algorithm}


\section{Error estimation}
\label{ch:error_est}
The adaptive MLSDC method requires a control mechanism for identifying regions where the spatial and/or temporal resolution needs to be adjusted. 
This task is performed by a set of error estimators that guide the refinement process. 
The total discretization error is assumed to consist of spatial, temporal, and hybrid contributions,
\begin{align}
    \varepsilon
    =
    \varepsilon_x
    +
    \varepsilon_t
    +
    \varepsilon_h.
\end{align}
Here, $\varepsilon_h$ denotes a hybrid error contribution arising from the interaction of spatial and temporal discretization errors. 
This contribution is not estimated explicitly. Numerical studies conducted by the authors based on analytical two-dimensional test functions indicate that $\varepsilon_h$ is negligible compared to the purely spatial and temporal error contributions. 
Furthermore, the proposed estimators are constructed to slightly overestimate the error, thereby providing a conservative refinement criterion.
A distinction should be made between \emph{error estimators}, which approximate the magnitude of the discretization error, and \emph{error indicators}, which merely identify regions requiring refinement without providing a quantitative error estimate. 
The estimators employed in this work and coefficient fitting procedures are derived from the works of Catherine Mavriplis \cite{EE_MavriplisDiss1989}.
Since the error estimation procedure is applied independently on every level $l$ and for every element $e$, the according indices are omitted throughout this section for improved readability.
It is assumed that a scalar sensor variable is available for error estimation. For scalar problems, this variable coincides with the solution itself, whereas for the Euler and compressible Navier-Stokes equations it corresponds to a selected solution component specified in the respective test case. For notational simplicity, the sensor variable is denoted by $u$ throughout this section.

For a given element $e$ and a single time step $n$, the nodal and modal representations of the space-time solution are related through the linear basis transformation
\begin{align}
\label{eq:lagrange_legendre}
\sum_{i=0}^{P}
\sum_{m=1}^{M}
{u}_{i,m}
\, \ell_i^x(\xi)\,\ell_m^t(\tau)
=
\sum_{i=0}^{P}
\sum_{m=1}^{M}
\hat{{u}}_{i,m}
\, P_i(\xi)\,P_m(\tau) .
\end{align}
Here, $P_i(\xi)$ and $P_m(\tau)$ denote the spatial and temporal Legendre basis functions, respectively, and $\hat{{u}}_{i,m}$ are the corresponding modal coefficients.

\subsection{Spatial Error}
The spatial error estimator almost exactly follows the approach proposed in \cite{EE_MavriplisDiss1989} and is applied to the spectral element approximation introduced in Section~\ref{sec:discretization}. 
The spatial error is evaluated exclusively at the final temporal collocation node $m=M$. Accordingly, the notation $\hat{u}_i \coloneqq \hat{u}_{i,M}$ is used throughout this spatial error subsection. Following \cite{EE_MavriplisDiss1989}, the estimator is composed of two contributions:

\medskip
\noindent\textbf{Truncation error.}
The first contribution estimates the truncation/ approximation error
\begin{align}
\|{u}-\breve{{u}}_h\|,
\end{align}
where $\breve{{u}}_h$ denotes the truncated polynomial expansion constructed from the exact modal coefficients. 
The notation here largely follows \cite{EE_MavriplisDiss1989}. 
This contribution measures the unresolved portion of the spectrum beyond the polynomial order $P$ (see Figure~\ref{fig:trunc_error_x}). Following \cite{EE_MavriplisDiss1989}, the truncation error may be approximated by
\begin{align}
    \|{u}-\breve{{u}}_h\|_{L^2}
    \approx
    \left(
        \sum_{i=P+1}^{\infty}
        \frac{2\hat{{u}}_i^2}{2i+1}
    \right)^{1/2}.
\end{align}
Since the modal coefficients $\hat{{u}}_i$ for $i>P$ are not explicitly available, their contribution must be approximated by a continuous representation. 
To this end, the discrete modal coefficients are assumed to admit a sufficiently smooth and slowly varying extension $\hat{{u}}(s)$ in the modal index such that $\hat{{u}}(s=i)\approx\hat{{u}}_i$ for integer values of $i$. 
The unresolved contribution is then approximated by replacing the discrete sum with the integral over the continuous mode variable $s$
\begin{align}
    \sum_{i=P+1}^{\infty}
    \frac{2\hat{ u}_i^2}{2i+1}
    \approx
    \int_{P+1}^{\infty}
    \frac{2\hat{ u}(s)^2}{2s+1}\,\mathrm{d}s .
\end{align}
For smooth solutions, the modal coefficients are assumed to exhibit exponential decay,
\begin{align}
    \hat{{u}}(s)
    \sim
    {c}\,e^{-\sigma s},
\end{align}
where the parameters ${c}$ and $\sigma$ are obtained from a least-squares fit using the last $N$ available Legendre coefficients of the current approximation. Unless stated otherwise, $N=4$ is used throughout this work.
For insufficiently smooth solutions, an algebraic decay model \cite{EE_MavriplisDiss1989}
may be employed instead. Both are implemented.
Note: the resulting truncation error estimate could also be used independently as an indicator for spatial $p$-adaptivity.

\medskip
\noindent\textbf{Quadrature/ interpolation error.}
The second contribution estimates the error
\begin{align}
\|{u}_h-\breve{{u}}_h\|,
\end{align}
which arises from the fact that the modal coefficients are only available through a numerical approximation. Following the terminology of \cite{EE_MavriplisDiss1989}, this contribution is referred to as the \emph{quadrature} or \emph{interpolation} error. 
In the spectral element setting, it additionally accounts for the non-optimality of the polynomial approximation, which is bounded by
\begin{align}
\|{u}_h-{u}^*\|_\infty
\le
\Lambda
\|{u}^*-{u}\|_\infty,
\end{align}
where ${u}^*$ denotes the best polynomial approximation and $\Lambda$ is the Lebesgue constant. 
A conservative estimate proposed by \cite{EE_MavriplisDiss1989} is given by
\begin{align}
    \|{u}_h-\breve{{u}}_h\|_{L^2}
    \approx
    \left(
        \frac{2\hat{{u}}_P^2}{2P+1}
    \right)^{1/2}.
\end{align}
This contribution may also be used independently as an indicator for spatial $h$-adaptivity.
Combining both contributions yields the spatial error estimate
\begin{align}
    {\varepsilon}_{\mathrm{est},x}
    =
    \left(
        \int_{P+1}^{\infty}
        \frac{2\hat{{u}}(s)^2}{2s+1}\,\mathrm{d}s
        +
        \frac{2\hat{{u}}_P^2}{2P+1}
    \right)^{1/2},
\end{align}
which is evaluated at the final temporal collocation node.

\subsection{Temporal Error}
The SDC method iteratively reduces the temporal discretization error by constructing a polynomial approximation in time at every spatial degree of freedom. 
Since this approximation is represented by Lagrange polynomials defined at the collocation nodes, it appears natural to consider transferring the spectral error estimator of \cite{EE_MavriplisDiss1989} used in space to the temporal discretization.
Upon closer inspection, however, the estimator is not directly applicable in the SDC context. The estimator of Mavriplis is designed to approximate the global $L^2$ error over an entire element interval. 
In contrast, the quantity of primary interest for SDC is the error at the final collocation node, where the collocation solution exhibits its characteristic superconvergence properties. 
Consequently, an estimator targeting the global temporal error does not necessarily provide an accurate measure of the error relevant for time-step advancement.

A second difficulty arises from the iterative nature of SDC. The modal coefficients obtained during intermediate sweeps do not correspond to an optimal polynomial projection. Instead, they evolve throughout the iterative process and only approach the collocation solution asymptotically as the number of sweeps increases.
For a finite number of sweeps, the solution at the collocation nodes typically exhibits an accuracy of approximately order $k$, where $k$ denotes the current sweep count. 
As a result, the quadrature and interpolation error contribution proposed by \cite{EE_MavriplisDiss1989}, which relies on projection-optimal modal coefficients, cannot be directly justified during the early stages of the SDC iteration.
Furthermore, the truncation error estimate of \cite{EE_MavriplisDiss1989} does not account for the superconvergence behavior at the final collocation node. For these reasons, a modified temporal error estimator is required for the adaptive MLSDC framework considered in this work.
Note: the temporal error is evaluated independently at each spatial node $i$. For notational simplicity, the spatial index is omitted throughout this subsection and the shorthand notation $\hat{u}_m \coloneqq \hat{u}_{i,m}$ is adopted. The aggregation of the resulting nodal error estimates into a element-wise temporal error contribution is discussed subsequently.

\medskip
\noindent\textbf{Extrapolation error.}
A different strategy is required for temporal error estimation for the SDC method. 
At low sweep counts, the quadrature/ interpolation error contribution is replaced by an estimator based on the difference between successive SDC approximations. 
This idea is closely related to Richardson extrapolation and is widely used in embedded Runge-Kutta methods, see, for example, \cite[Chap.~II.4]{TI_Hairer1993a}. A similar approach has recently been employed for time-step adaptivity in SDC methods by \cite{EE_Baumann2024}.
We refer to this contribution as the \emph{extrapolation error} ${\varepsilon}_{\mathrm{ext},t}$, which measures the discrepancy between the two most recent approximations at the final collocation node,
\begin{align}
    {\varepsilon}_{\mathrm{ext},t}
    =
    \left\|
        {u}_{M}^{k_{\max}}
        -
        {u}_{M}^{k_{\max}-1}
    \right\|_{L^2}.
\end{align}
Here, the $L^2$ norm is evaluated over the spatial domain.
The extrapolation error dominates during the early stages of the SDC iteration, when the polynomial coefficients are still evolving toward the collocation solution. 
However, this contribution decreases monotonically with increasing sweep count and vanishes as the SDC iteration converges. 
Consequently, it cannot serve as a complete temporal error indicator on its own, since it eventually becomes insensitive to the remaining temporal discretization error.

\begin{figure}[H]
    \centering
    \begin{subfigure}{\textwidth}
        \centering
        \includegraphics[width=0.8\textwidth]{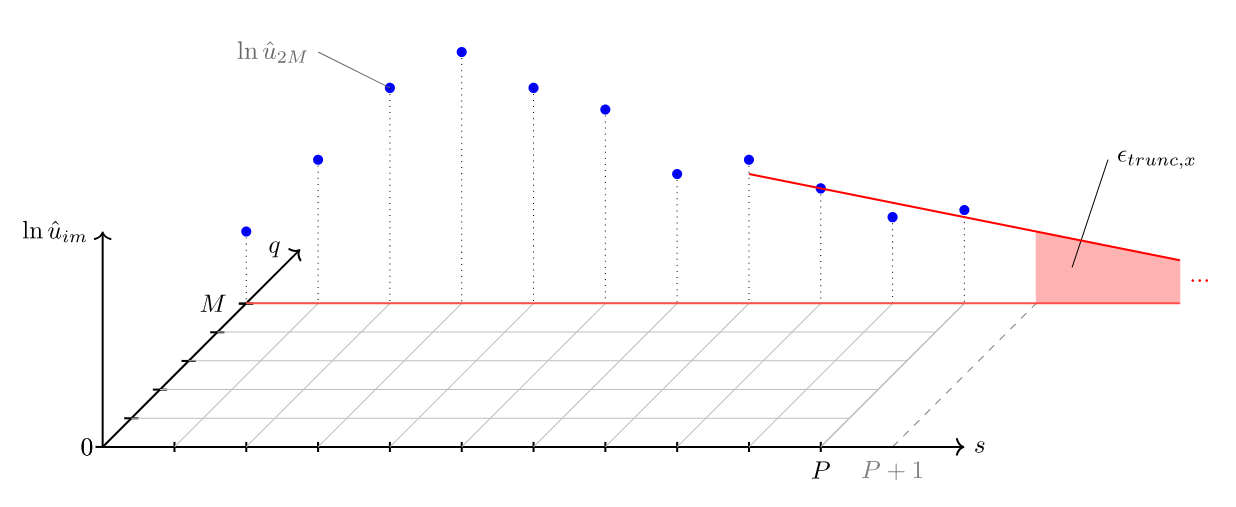}
        \caption{Spatial truncation error estimate at temporal node $m=M$.}
        \label{fig:trunc_error_x}
    \end{subfigure}

    \vspace{0.3cm}

    \begin{subfigure}{\textwidth}
        \centering
        \includegraphics[width=0.8\textwidth]{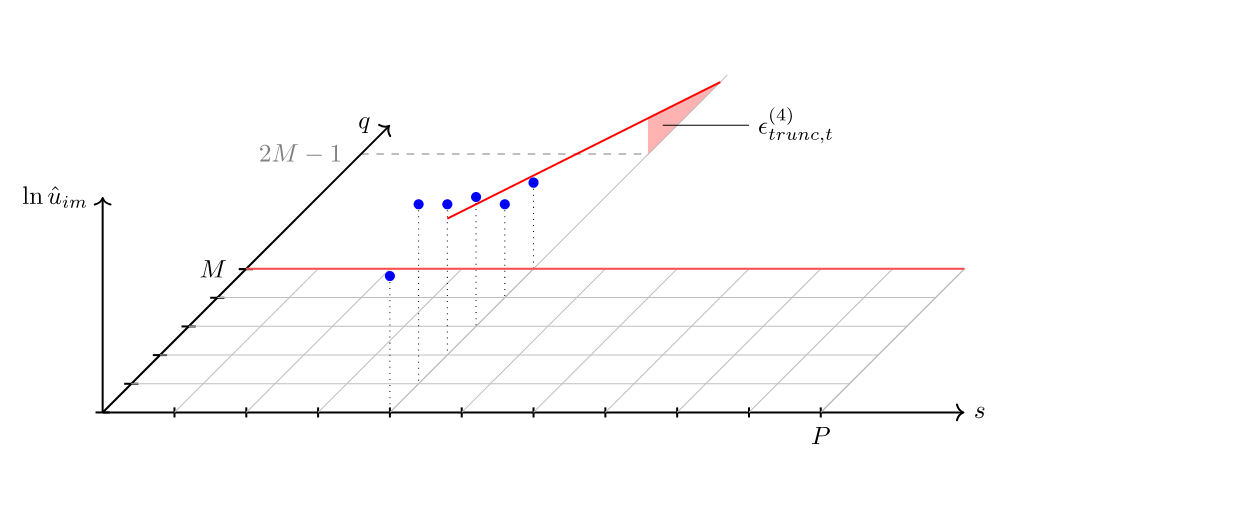}
        \caption{$i$-th temporal truncation error estimate.}
        \label{fig:trunc_error_t}
    \end{subfigure}
    \caption{Sketches of the calculation of the truncation error estimates.}
    \label{fig:trunc_error}
\end{figure}

\medskip
\noindent\textbf{Modified truncation error.}
For sufficiently many sweeps, the temporal discretization error is instead dominated by the unresolved temporal spectrum. 
In this regime, the truncation error is estimated analogously to the approach of \cite{EE_MavriplisDiss1989}. However, the superconvergence of the final collocation node is incorporated through a small modification.
Rather than estimating the unresolved spectrum from mode $M+1$ onward, the integration is performed over the interval $[2M-1,\infty)$ (see Figure~\ref{fig:trunc_error_t}). 
This modification is motivated by the fact that SDC with Radau-Right nodes converges at the last collocation node to the Radau IIA collocation method, whose order of accuracy is $2M-1$ \cite[Ch.~6]{TI_Deuflhard2002a} (superconvergence).
While the $M$ collocation nodes determine the coefficients within the polynomial approximation space exactly, the superconvergent Radau endpoint achieves an accuracy of order $2M-1$. Consequently, the leading unresolved contribution to the error is expected to be associated with modes beyond order $2M-1$. 
We refer to the resulting contribution as the \emph{modified truncation error} ${\varepsilon}_{\mathrm{trunc},t}$.
As in the spatial case, the discrete temporal modal coefficients $\hat{{u}}_m$ are assumed to admit a sufficiently smooth extension $\hat{{u}}(q)$ with respect to the temporal mode index $q$. The fitting procedure used to obtain $\hat{{u}}(q)$ is identical to that employed for the spatial estimator. The modified truncation error is then defined by
\begin{align}
    {\varepsilon}_{\mathrm{trunc},t}
    =
    \int_{2M-1}^{\infty}
    \frac{2 \hat{{u}}(q)^2}{2q+1}\,\mathrm{d}q .
\end{align}
The modified truncation error becomes the dominant contribution as the number of sweeps increase and provides a reliable lower bound for the remaining temporal discretization error. 
Moreover, this contribution may be employed independently as an indicator for temporal $m$-adaptivity.

The overall temporal error estimate is obtained by taking the maximum of the the extrapolation and modified truncation errors,
\begin{align}
    {\varepsilon}_{\mathrm{est},t}
    =
    \max\!\left(
        {\varepsilon}_{\mathrm{trunc},t},
        {\varepsilon}_{\mathrm{ext},t}
    \right).
\end{align}
This temporal estimate is evaluated at each spatial node $i$ and then the $L^2$ norm over all spatial degrees of freedom within an element is then computed.
The overall error estimate is then defined as
\begin{align}
    \varepsilon_{\mathrm{est}}
    =
    \varepsilon_{\mathrm{est},x}
    +
    \varepsilon_{\mathrm{est},t}.
\end{align}
The spatial and temporal contributions are combined additively in order to obtain a conservative estimate of the total discretization error.
\begin{remark}
The error estimates obtained using the approach suggested by \cite{EE_MavriplisDiss1989} are computed on the standard element. To make the spatial or temporal contribution directly comparable to the exact element-wise $L_2$ error introduced below, the corresponding metric factor is incorporated.
\end{remark}
\begin{remark}
The accuracy of the proposed spectral error estimators depends critically on the quality of the modal coefficients used for the least-squares fitting procedure. 
In practice, numerical noise, round-off errors, and weakly resolved solution features may contaminate the highest modal coefficients and lead to unreliable decay-rate estimates. Consequently, a filtering procedure is applied prior to the coefficient fit.
A simple filtering strategy removes coefficients whose magnitude falls below a prescribed threshold. In addition, an advanced filtering procedure is employed to suppress oscillatory coefficient patterns that are inconsistent with the assumed smooth modal decay. Such oscillations typically arise from numerical noise and may significantly deteriorate the quality of the fitted decay model.
If an excessive number of coefficients must be filtered, the remaining information is insufficient for a reliable fit. In this case, the estimator automatically falls back to the approach proposed by \cite{EE_Henderson1999}.
\end{remark}

\subsection{Test}
\begin{figure}[H]
    \centering
    \begin{subfigure}[t]{0.48\textwidth}
        \centering
        \includegraphics[width=\textwidth]{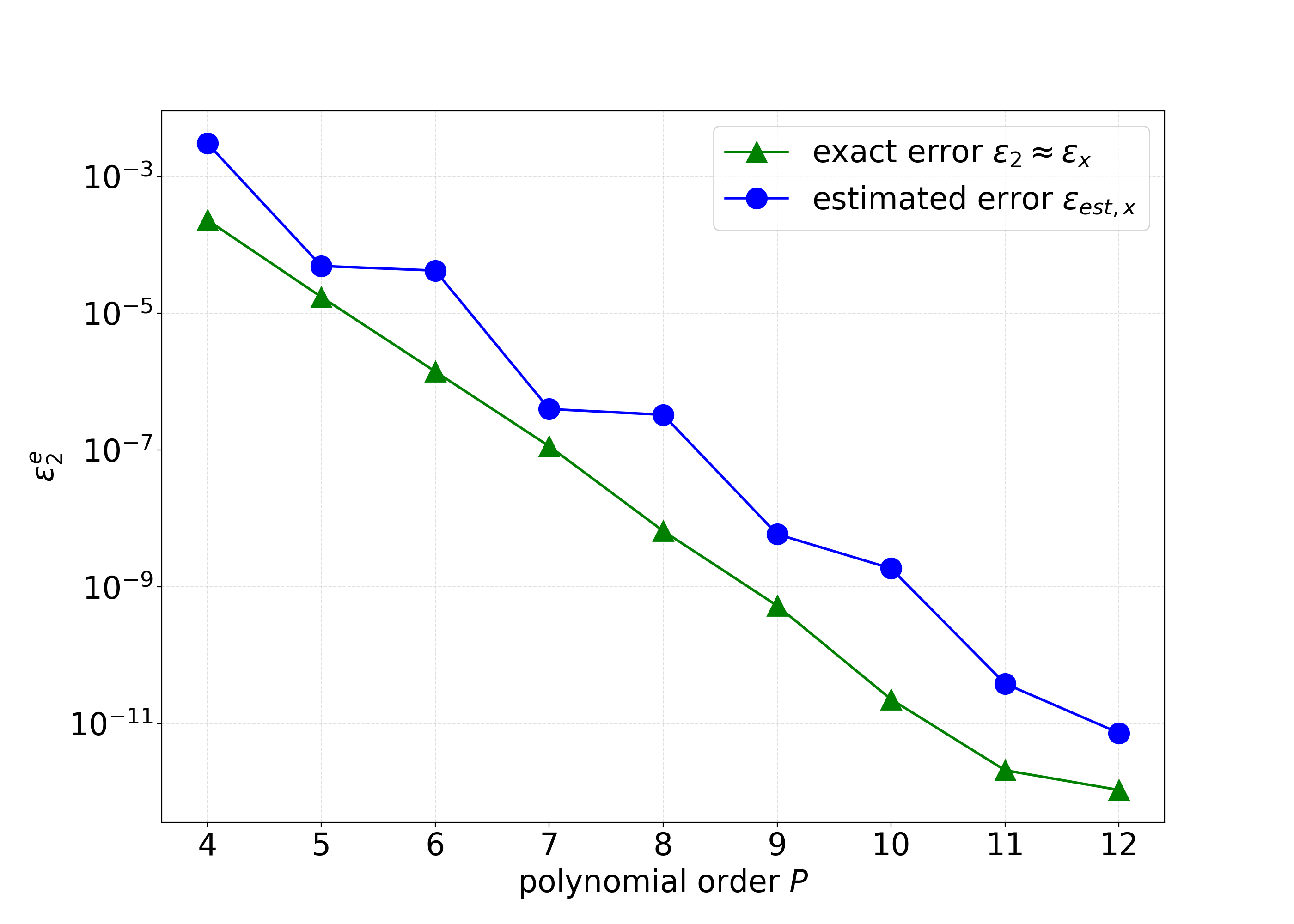}
        \caption{Exact vs estimated error for spatially \\ dominated problem.}
        \label{fig:errors_est_x}
    \end{subfigure}%
    \begin{subfigure}[t]{0.48\textwidth}
        \centering
        \includegraphics[width=\textwidth]{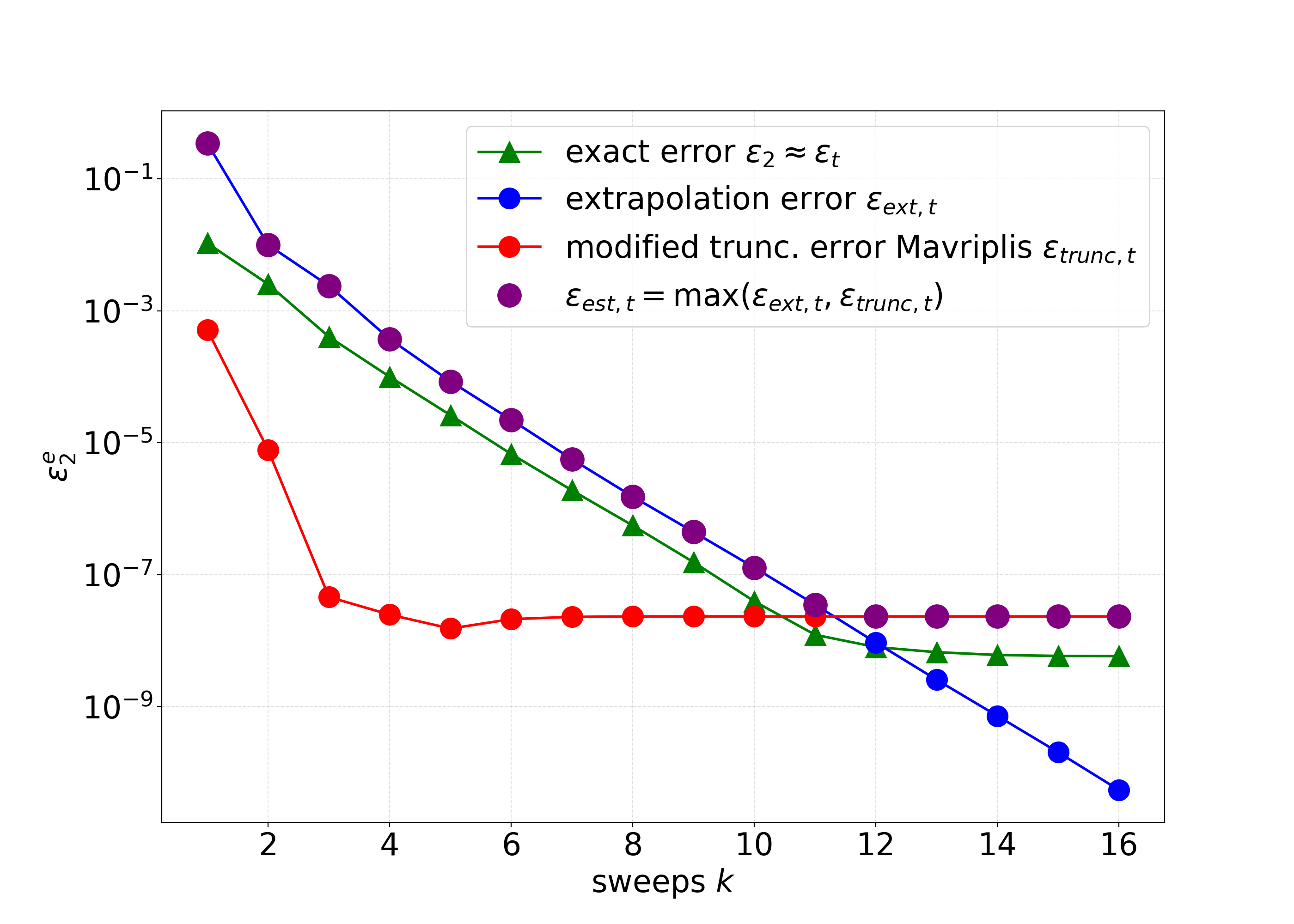}
        \caption{Exact vs extrapolation vs spectral error for temporally dominated problem.}
        \label{fig:errors_est_t}
    \end{subfigure}
    
    \caption{Performance of spatial and temporal error estimators for the convection-diffusion wave package problem.}
    \label{fig:err_combined}
\end{figure}

The performance of the combined spatial and temporal error estimator is assessed in Figure~\ref{fig:err_combined} using the convection-diffusion equation \eqref{eq:conv-diff} with diffusivity $\nu=0.02$ on the domain $\Omega=[0,1]$ for single-level SDC. The analysis is performed over a single time step. The exact solution is given by
\begin{equation}
  \label{eq:wave-packet}
  u(x,t) = \sum_{i=1}^7 a_i \sin\big( \kappa_i(x - s_i - vt) \big)
                        e^{-\kappa_i^2 \nu t}
  .
\end{equation}
Equation~\eqref{eq:wave-packet} defines a wave packet consisting of multiple wave numbers $\kappa_i$, amplitudes $a_i$, and phase shifts $s_i$, whose values are listed in Table~\ref{tab:wave-package}. All simulations employ the SI(1) time integrator within the SDC method.
\begin{table}[H]
  \begin{tabular}{lccccccc} 
  \toprule
  $i$        & $1$    & $2$    & $3$     & $4$     & $5$     & $6$     & $7$     \\\midrule
  $\kappa_i$ & $2\pi$ & $6\pi$ & $10\pi$ & $14\pi$ & $18\pi$ & $24\pi$ & $30\pi$ \\
  $a_i$      & $1.00$ & $1.50$ & $1.80$  & $1.70$  & $1.50$  & $1.30$  & $1.15$  \\
  $s_i$      & $0.00$ & $0.05$ & $0.10$  & $0.15$  & $0.20$  & $0.30$  & $0.18$  \\
  \bottomrule
  \end{tabular}
    \caption{Wave package coefficients.}
    \label{tab:wave-package}
\end{table}
The exact element-wise $L^2$ error is measured by
\begin{align}
    \epsilon^{e}_2
    =
    \left\|
        u(x,t^{n})
        -
        \NM {u}_M^e(x,t^n)
    \right\|_{L^2},
\end{align}
where $\NM {u}_M^e$ denotes the numerical solution at the final collocation node within element $\Omega^e$.
Figure~\ref{fig:errors_est_x} considers a spatially dominated error regime. To minimize the influence of temporal discretization errors, a time-step size of $\Delta t=10^{-3}$ and $12$ SDC sweeps are employed. 
The other discretization parameters can be taken from Table~\ref{tab:discretization_params}.
\begin{table}[H]
\centering
\begin{tabular}{ccccc}
\toprule
 & \multicolumn{1}{c}{Space} & \multicolumn{2}{c}{Time} \\
\cmidrule(r){2-2} \cmidrule(r){3-4}
 & $N_{e,l}$ & $n_{t,l}$ & $M_l$ \\
\midrule
Level 1 & 32 & 1 & 6 \\
\bottomrule
\end{tabular}
\caption{Discretization parameters used both tests evaluating the performance of the error estimator.}
\label{tab:discretization_params}
\end{table}
Under these conditions, the measured error is governed primarily by the spatial discretization, allowing an isolated assessment of the spatial error estimator.
Figure~\ref{fig:errors_est_x} shows the estimated and exact errors for one representative element of the computational domain. 
Similar results were observed for all other remaining elements. 
Results are presented for a range of polynomial orders. The estimator provides a consistent overestimate of the exact error while accurately reproducing its qualitative behavior. The observed plateaus originate from symmetries in the wave-packet solution and are therefore not directely related to deficiencies of the estimator. \\
The temporal error estimator is assessed in Figure~\ref{fig:errors_est_t}. To obtain a temporally dominated error regime, the simulations are performed with $\Delta t=10^{-2}$ and a fixed polynomial degree of $P=15$. 
The other discretization parameters can be taken from Table~\ref{tab:discretization_params}.
Again, the results shown correspond to a representative element. Good agreement between the estimated and exact errors is observed throughout the SDC iterations. 
During the initial sweeps, the extrapolation error dominates, reflecting the ongoing convergence of the SDC solution toward the collocation solution. Note on how the SDC method converges after $11$ sweeps with $M=6$.
For larger sweep counts, the modified truncation error becomes the dominant contribution. The temporal error estimate, defined as the maximum of these two contributions, captures the exact error remarkably well over the entire range of sweep counts.
The results shown in Figure~\ref{fig:errors_est_t} are obtained using a single-level SDC scheme. However, the error estimation strategy is equally applicable to MLSDC methods, since the estimator is evaluated independently on each level. The only practical difference is that the modal coefficients generally converge more rapidly toward the collocation solution in the multilevel setting.
\begin{remark}
    For single-level SDC, the modified temporal truncation error no longer changes once $k=M$, since the errors at the collocation nodes stagnate beyond this point. 
\end{remark}


\section{Numerical Studies}
\label{ch:tests}

\subsection{Burgers moving front}
\begin{figure}[H]
    \centering
    \resizebox{0.9\textwidth}{!}{%
    \begin{minipage}{\textwidth}
    \centering   
    \begin{subfigure}{\linewidth}
     \centering
     \includegraphics[width=\linewidth]{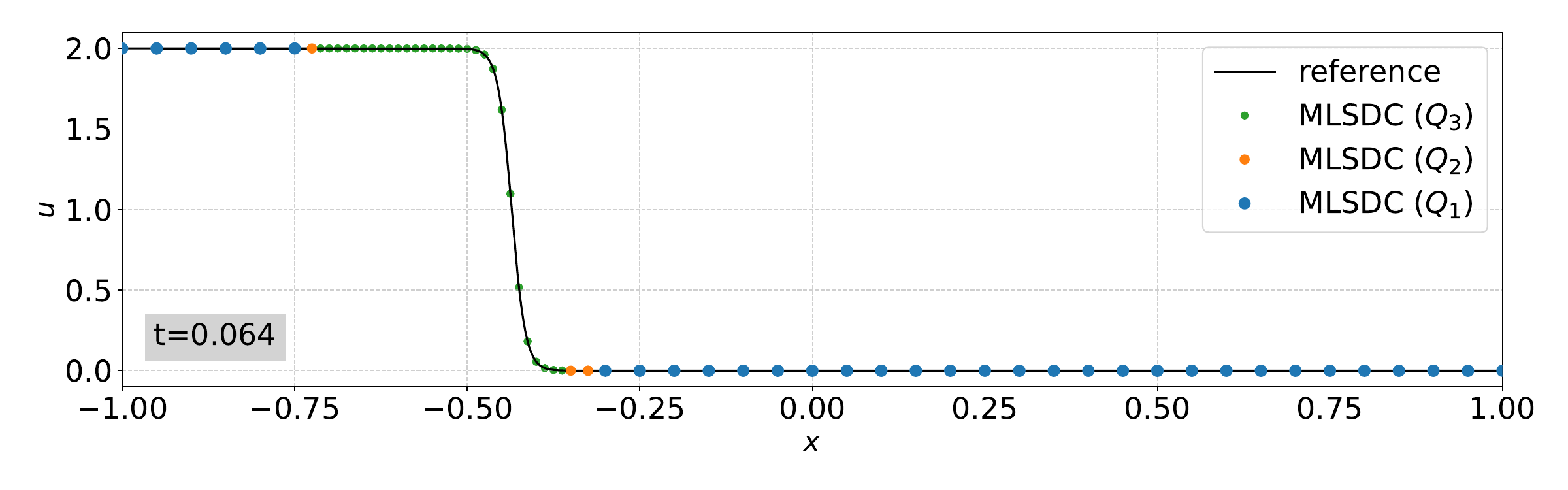}
     \caption{Solution at $t = 0.064$.}
    \end{subfigure}
    \vspace{0.5em}
    \begin{subfigure}{\linewidth}
        \centering
        \includegraphics[width=\linewidth]{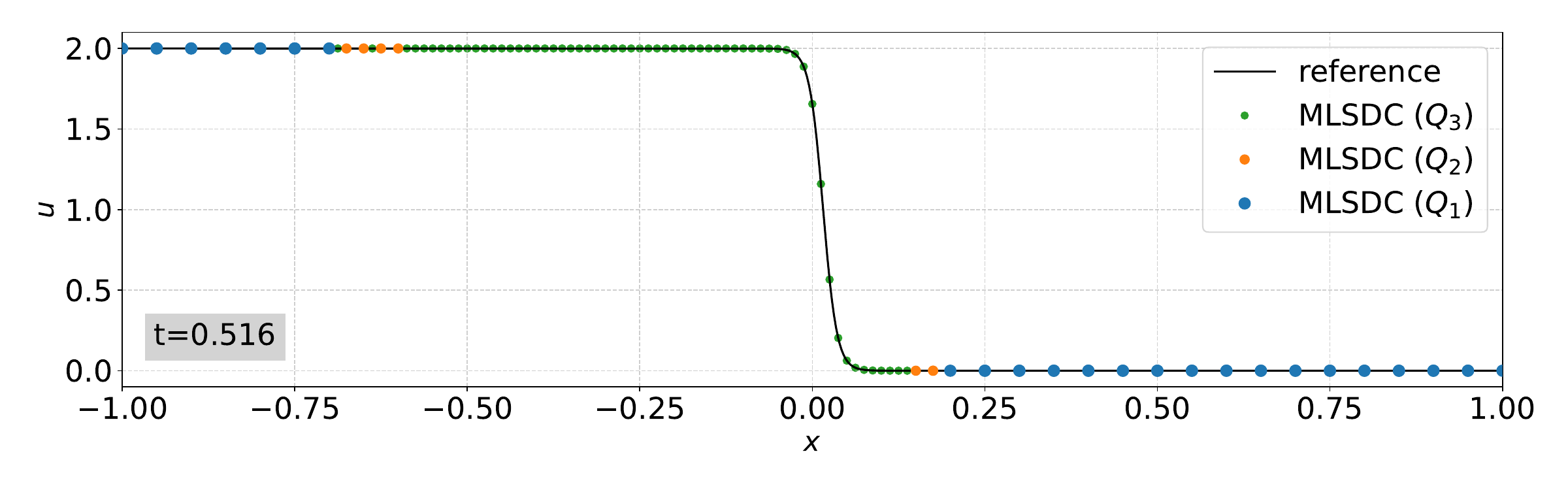}
        \caption{Solution at $t = 0.516$.}
    \end{subfigure}
    \vspace{0.5em}
    \begin{subfigure}{\linewidth}
        \centering
        \includegraphics[width=\linewidth]{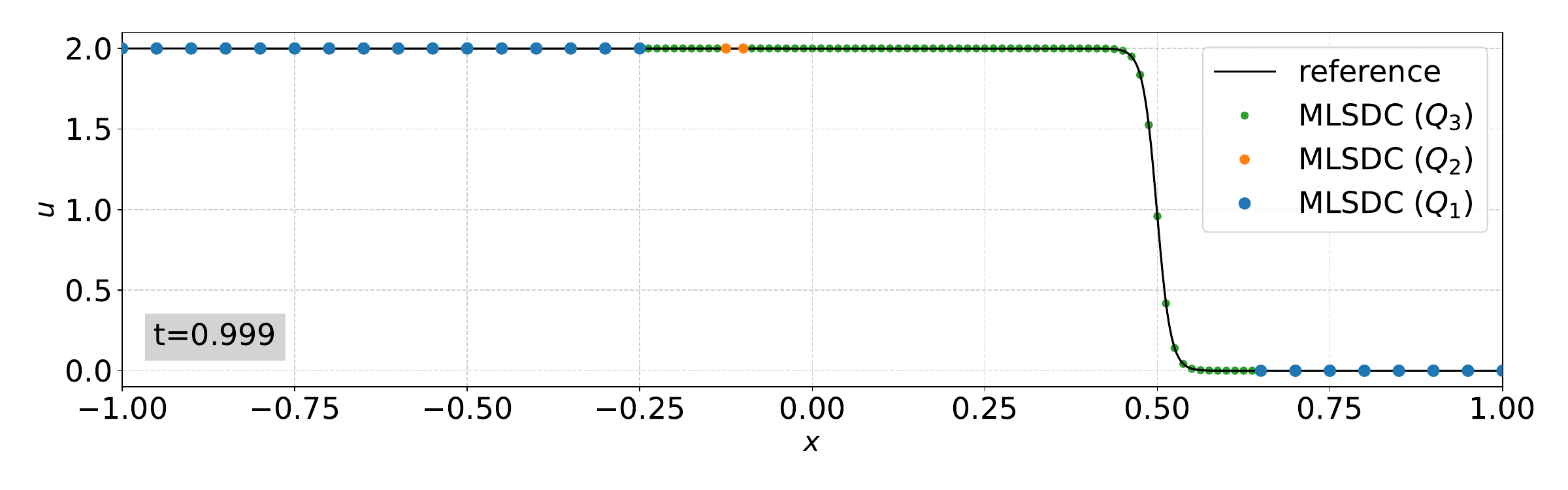}
        \caption{Solution at $t = 0.999$.}
     \end{subfigure}
    \end{minipage}%
    }
    \caption{Burgers moving front solved with adaptive MLSDC$_{3,5,7}^2$-SI(1) at different time steps.}
    \label{fig:burgers_MF_progression}
\end{figure}

\begin{figure}[H]
    \centering
    \includegraphics[width=0.8\textwidth]{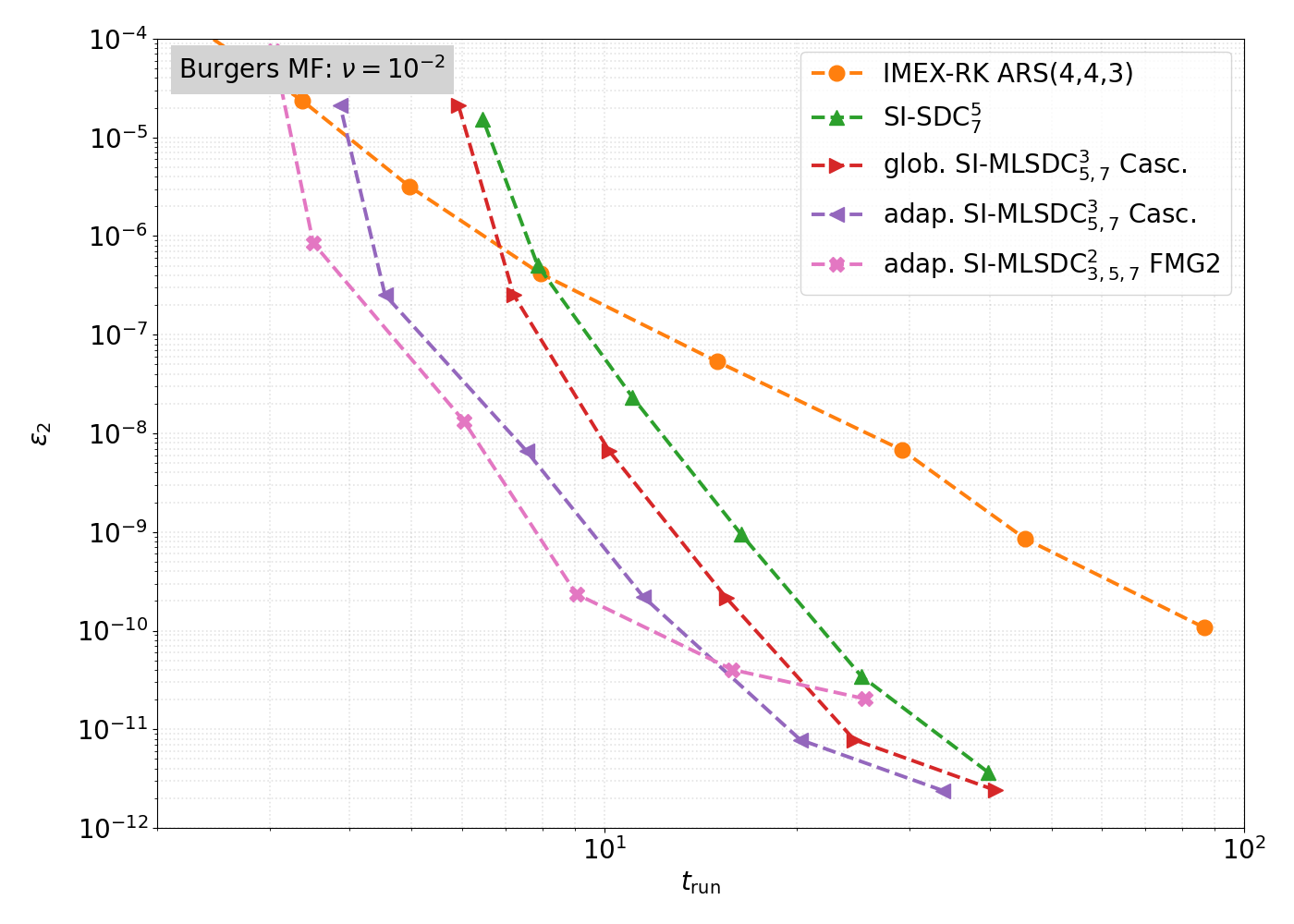}
    \caption{Error-runtime comparison between a Runge-Kutta, SDC and locally as well as globally refined MLSDC methods using either the preconditioned conjugate gradient method (PCG) or the fast direct solver (FDS) for the Burgers moving front problem with $\nu = 10^{-2}$. }
    \label{fig:burgers_mf_exact}
\end{figure}

The capabilities of the adaptive MLSDC$_{3,5,7}^{2}$ method are assessed using the proposed error estimator. Here, the notation MLSDC$_{3,5,7}^{2}$ denotes a three-level hierarchy employing $3$, $5$, and $7$ temporal collocation nodes on the respective levels and performing two multilevel V-cycles per time step. The FMG starting stragety was used for initialization and there was no post-sweeping step used after the last V-cycle.
As a first test case, the Burgers equation \eqref{eq:burgers} is solved. The computational domain is given by $\Omega=[-1,1]$, the viscosity is set to $\nu=0.01$, and the simulation is carried out until $t_{\mathrm{end}}\approx1$ using a time-step size of $\Delta t\approx0.032$. The spatial and temporal discretization parameters are summarized in Table~\ref{tab:num_experiments_burgers_runtime}.
    \begin{table}[H]
    \centering
    \begin{tabular}{cccc}
    \toprule
     & \multicolumn{2}{c}{Space} & \multicolumn{1}{c}{Time} \\
    \cmidrule(r){2-3} \cmidrule(r){4-4}
     & $N_{e,l}$ & $P_l$ & $M_l$ \\
    \midrule
    Level 1 & 40 & 12 & 3 \\
    Level 2 & 80 & 12 & 5 \\
    Level 3 & 160 & 12 & 7 \\
    \bottomrule
    \end{tabular}
    \caption{Discretization parameters for the Burgers runtime tests.}
    \label{tab:num_experiments_burgers_runtime}
    \end{table}
\begin{remark}
The authors are aware that the choice of coordinated space-time refinement and prescribed finer-level discretizations may appear somewhat arbitrary. This strategy was adopted deliberately in order to isolate and investigate the effects of locally increased space-time resolution.
\end{remark}
The novel error estimator from Section~\ref{ch:error_est} serves as reference for error evaluation. Adaptive refinement and coarsening are triggered using thresholds of $10^{-8}$ and $10^{-11}$, respectively.
The solution exhibits a steep front propagating from left to right throughout the simulated time. Figure~\ref{fig:burgers_MF_progression} illustrates the corresponding adaptive refinement pattern at $t\approx0$, $t\approx0.5$, and $t\approx1$, demonstrating the ability of the error estimator to accurately track the moving front.
All simulations are performed using the SI(1) time integrator for predictor and corrector within the MLSDC$_{3,5,7}^{2}$ method.
The points shown in Figure~\ref{fig:burgers_MF_progression} indicate element boundaries and therefore visualize the dynamically evolving mesh hierarchy.
The results demonstrate that the moving front is detected reliably by the error estimator. As the front propagates through the domain, the algorithm adapts the mesh accordingly and coarsens previously refined regions once high resolution is no longer required.
A slight delay in the coarsening process can be observed in the wake of the refinement zone. This behavior may be attributed to a hysteresis-like effect caused by the separation of the refinement and coarsening thresholds.
Due to the pronounced separation between regions of high and low error, the intermediate refinement level is used only sparsely in this example. The solution is effectively characterized by a small region of large error and a large region of negligible error, resulting in a predominantly two-level behavior.
Compared to a uniformly refined single-level discretization achieving the same accuracy, the adaptive MLSDC method yields a substantial reduction in computational cost, requiring only approximately $60\%$ of the runtime.
The runtime savings obtained through adaptive refinement are quantified more systematically in the following experiment.
To this end, a performance study is conducted comparing error and runtime characteristics of several MLSDC methods as well as conventional time-integration methods.
Figure~\ref{fig:burgers_mf_exact} presents runtime-to-solution diagrams comparing IMEX RK-ARS(4,4,3) \cite{TI_Ascher1997a}, single-level SDC, MLSDC, and adaptive MLSDC. 
The runtime experiments were performed with the HiSPEET framework compiled using the Intel Fortran Compiler 2021.10.0, build 20230609. All runs were executed in single-threaded mode without parallelization on an Intel Xeon Platinum 8470 CPU  operating at 2.00\,GHz. Each configuration was repeated ten times, and the reported runtime corresponds to the arithmetic mean over these runs.
The results show a clear performance improvement when moving from classical Runge-Kutta methods to single-level SDC, globally refined MLSDC, and finally locally adaptive MLSDC.
In this test case, the TVD-RK3 method \cite{TI_Shu1988b} is not competitive with the SDC-based approaches, confirming the results shown in Section~8.4 in \cite{MG_Pfister2025}.
At very small error levels, mild deviations from the expected runtime-to-accuracy behavior can be observed, particularly for the three-level MLSDC method. These deviations are likely caused by the accumulation of round-off errors due to the increased number of algebraic operations in the multilevel iteration.
The current MLSDC implementation has not been optimized specifically for runtime. Further improvements are expected through parameter tuning, for example of the refinement and coarsening thresholds, since the performance of adaptive MLSDC methods depends on several algorithmic parameters.
Nevertheless, the results already demonstrate substantial potential for locally adaptive MLSDC methods and provide a strong motivation for future extensions to three-dimensional problems.
The observed speedup is caused by two complementary effects: the multilevel correction strategy already accelerates the globally refined MLSDC method relative to single-level SDC, while local adaptivity further reduces the computational work by concentrating fine-level resolution only where it is required.

\subsection{Refinement criteria}
\label{sec:refine}
The moving-front Burgers test case introduced in the previous subsections is reconsidered with a MLSDC$_{3,5,7}^2$.
The objective of this study is to compare different refinement criteria. Up to this point, only the novel proposed space-time error estimator has been employed to drive adaptive refinement.
To this end, the exact $L^2$ error, the proposed error estimator, the estimator of Henderson \cite{EE_Henderson1999}, and a refinement criterion based on the artificial diffusivity sensor of Persson and Peraire \cite{SE_Persson2006a} are considered. 
For all experiments, the three-level MLSDC method is used together with the discretization parameters listed in Table~\ref{tab:num_experiments_burgers_runtime}. Runtime and achieved $L^2$ error serve as the primary performance metrics.
For the exact error, the proposed error estimator, and the Henderson estimator, an element is refined whenever the corresponding error indicator exceeds a prescribed refinement threshold. For the artificial-diffusivity criterion, refinement is triggered whenever a nonzero artificial diffusivity is detected.
Coarsening is performed whenever the estimated error falls below a prescribed coarsening threshold. Assuming that the error is dominated by the spatial discretization and considering pure $h$-refinement, the relationship between coarse- and fine-grid errors may be approximated by
\begin{align}
\varepsilon_c \approx 2^{P+1}\varepsilon_f,
\end{align}
where $P$ denotes the spatial polynomial degree.
For $P=12$, this estimate predicts a coarsening factor of $2^{13}=8192$. However, numerical experiments indicate that a smaller factor of approximately $1000$ yields superior performance for the present test case.
It should be noted that the adaptive algorithm contains several tunable parameters whose further optimization may lead to additional performance improvements.
All runtime experiments were performed using the HiSPEET framework compiled with the Intel Fortran Compiler 2021.10.0 (build date 20230609). Simulations were executed in single-threaded mode on an Intel Xeon Platinum 8470 CPU operating at 2.00,GHz. Each configuration was repeated ten times, and the reported runtime corresponds to the average over all runs.
\begin{table}[H]
\centering
\begin{tabular}{c c c c c}
\hline
refine & coarsen & runtime & runtime rel. & error rel. \\
\hline
$10^{-14}$ & $10^{-17}$ & 6.026 & 0.692 & 1.000 \\
$10^{-13}$ & $10^{-16}$ & 5.793 & 0.665 & 1.000 \\
$10^{-12}$ & $10^{-15}$ & 5.533 & 0.635 & 1.001 \\
$10^{-11}$ & $10^{-14}$ & 5.123 & 0.588 & 1.014 \\
$10^{-10}$ & $10^{-13}$ & 4.078 & 0.468 & 2.619 \\
$10^{-9}$ & $10^{-12}$ & 3.914 & 0.449 & 5.575 \\
$10^{-8}$ & $10^{-11}$ & 3.739 & 0.429 & 135.440 \\
\hline
\end{tabular}
\caption{Full threshold table for the criterion exact error. The time slab size is $\Delta t = 8.058 \times 10^{-3}$. Refinement thresholds range from $10^{-14}$ to $10^{-8}$ in powers of ten. Runtime is given relative to the global-refinement runtime (8.710\,s), and error relative to the best case within this file.}
\label{tab:exact_error_thresholds}
\end{table}
The number of elements marked for refinement depends on both the chosen refinement criterion and the corresponding refinement threshold. Consequently, both the achieved accuracy and the computational cost are directly influenced by these parameters.
In addition there is theoretically a computational overhead by the estimators/indicators. For the methods under consideration however, they are negligible.
For each refinement criterion, a series of simulations is performed over a range of refinement thresholds together with the associated coarsening thresholds. This procedure is repeated for several time-step sizes in order to assess the robustness and performance of the refinement strategy.
Table~\ref{tab:exact_error_thresholds} presents a representative example based on the exact error criterion.
To facilitate a direct comparison of the considered refinement criteria, Table~\ref{tab:summary_thresholds_dt0.008} summarizes the refinement and coarsening threshold combinations yielding the best performance under the constraint that the resulting error does not exceed $1.1$ times the error obtained on the globally refined mesh.
The dominant differences in runtime originate from the number of elements selected for refinement and the resulting distribution of computational work.
It should furthermore be noted that the artificial-diffusivity criterion acts solely as an error indicator, whereas the remaining approaches provide genuine error estimates that can be directly related to the discretization error.
\begin{table}[H]
\centering
\begin{tabular}{l c c c c c}
\hline
criterion & refine & coarsen & runtime & runtime rel. & error rel. \\
\hline
Exact error & $10^{-11}$ & $10^{-14}$ & 5.123 & 0.588 & 1.014 \\
Error estimator & $10^{-10}$ & $10^{-13}$ & 4.178 & 0.480 & 1.002 \\
Henderson & $10^{-10}$ & $10^{-13}$ & 4.126 & 0.474 & 1.003 \\
Artificial diffusivity & - & - & 8.460 & 0.971 & 1.000 \\
\hline
\end{tabular}
\caption{Fastest acceptable configuration for each refinement criterion with relative error $\le 1.10$ relative to the error for global refinement. The slab size is $\Delta t = 8.058 \times 10^{-3}$. Refinement thresholds range from $10^{-14}$ to $10^{-8}$ in powers of ten. Runtime is given relative to the global-refinement runtime (8.710\,s).}
\label{tab:summary_thresholds_dt0.008}
\end{table}
The artificial-diffusivity criterion consistently overestimates the extent of the refinement zone. While this behavior may, in principle, be mitigated through parameter tuning, preliminary studies with different scaling factors showed a clear trade-off. The detected refinement region either remained excessively large, reducing efficiency, or became too narrow to adequately resolve the wave front, resulting in significantly increased errors.
Furthermore, any choice of these parameters would be largely ad hoc, given that the artificial diffusivity constitutes only a heuristic error indicator. The criterion remains also fundamentally limited to the detection of spatial under-resolution and is therefore expected to lose effectiveness once temporal errors become significant.
Both the proposed error estimator and the exact error criterion identify compact refinement regions that closely track the relevant solution features and therefore yield efficient adaptive meshes. They show significant runtime reductions while maintaining accuracy for this test case.
Obviously the exact error is generally unavailable in practical applications and therefore cannot normally be used as an adaptive refinement criterion. It only serves as a comparison here.
Among the considered approaches, the estimator of Henderson \cite{EE_Henderson1999} achieves the lowest runtime while maintaining an acceptable error level. However, its applicability is problem dependent, and limitations become apparent in more challenging test cases, see Section~\ref{subsec:shu_osher}.
The observed trends are reproduced for a smaller time-step size, as demonstrated by the results summarized in Table~\ref{tab:summary_thresholds_dt0.004}.
\begin{table}[H]
\centering
\begin{tabular}{l c c c c c}
\hline
criterion & refine & coarsen & runtime & runtime rel. & error rel. \\
\hline
Exact error & $10^{-12}$ & $10^{-15}$ & 8.204 & 0.626 & 1.037 \\
Error estimator & $10^{-11}$ & $10^{-14}$ & 6.460 & 0.493 & 1.046 \\
Henderson & $10^{-10}$ & $10^{-13}$ & 6.200 & 0.473 & 1.007 \\
Artificial diffusivity & - & - & 12.138 & 0.926 & 1.000 \\
\hline
\end{tabular}
\caption{Fastest acceptable configuration for each refinement criterion with relative error $\le 1.10$ relative to the error for global refinement. The slab size is $\Delta t = 4.029 \times 10^{-3}$. Refinement thresholds range from $10^{-15}$ to $10^{-9}$ in powers of ten. Runtime is given relative to the global-refinement runtime (13.110\,s).}
\label{tab:summary_thresholds_dt0.004}
\end{table}

\subsection{Shu-Osher shock fluctuation benchmark}
\label{subsec:shu_osher}
\begin{figure}[H]
    \centering
    \begin{subfigure}{0.48\textwidth}
        \centering
        \includegraphics[width=\textwidth]{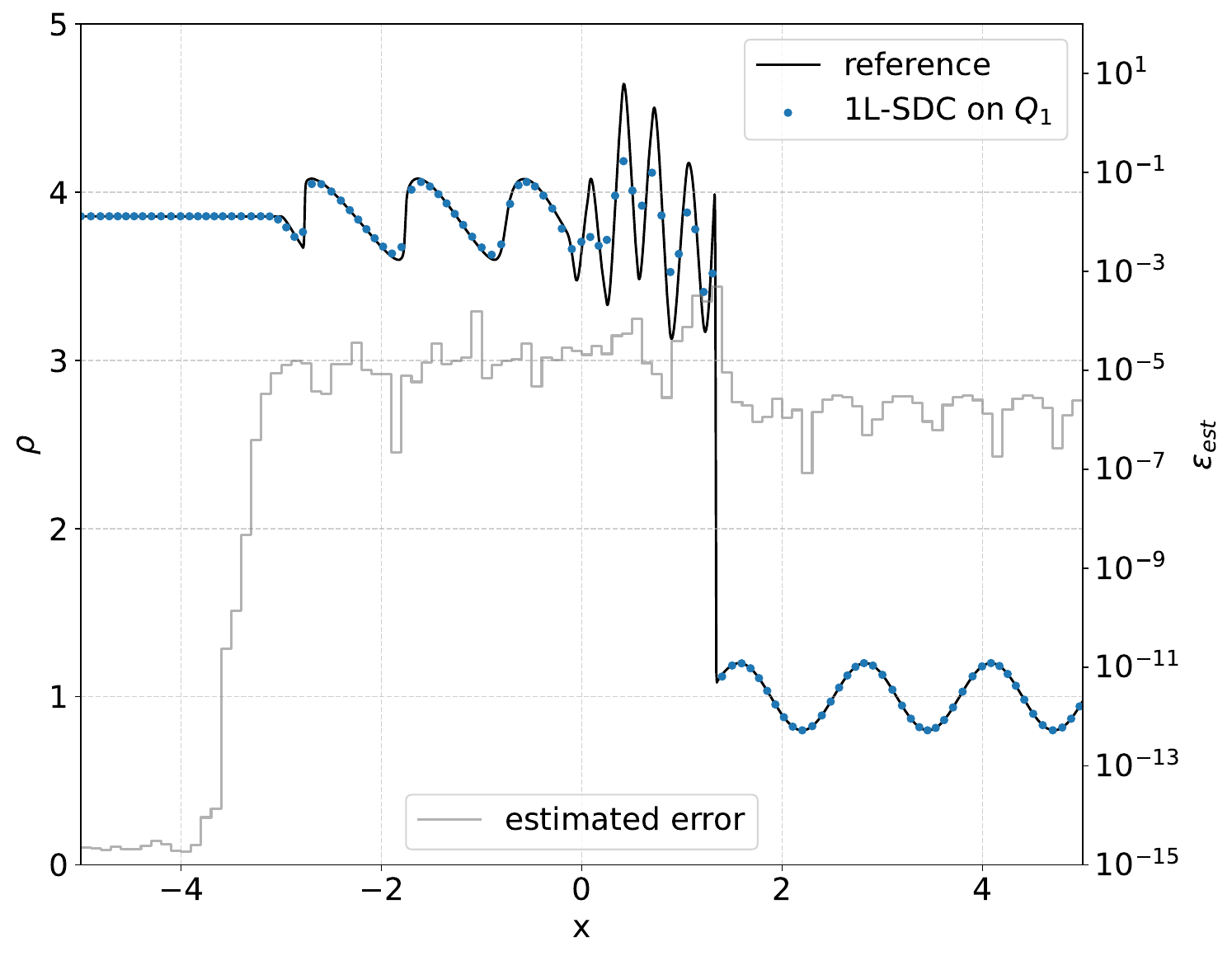}
        \caption{Solution for single-level SDC (level 1).}
    \end{subfigure}
    \hfill
    \begin{subfigure}{0.48\textwidth}
        \centering
        \includegraphics[width=\textwidth]{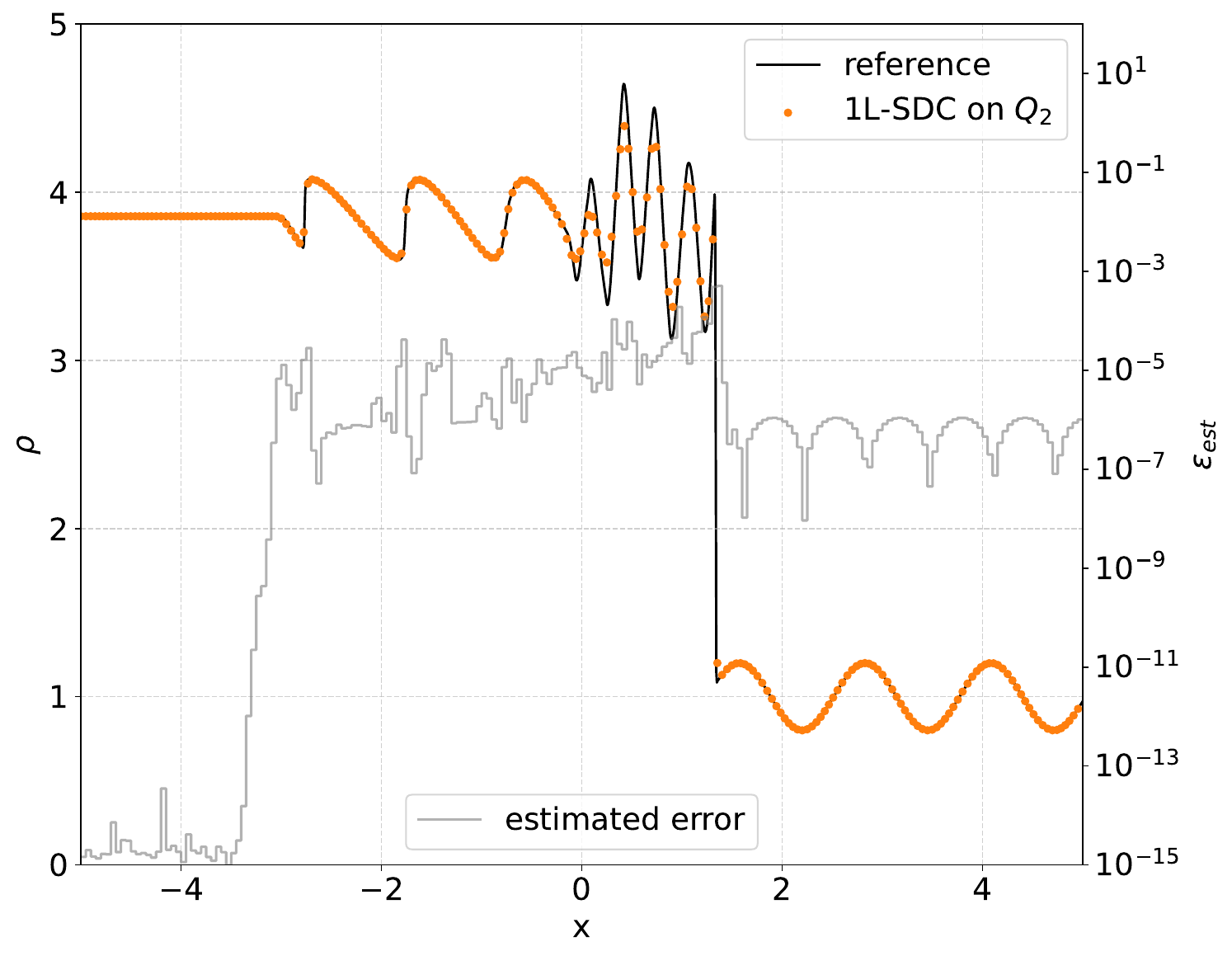}
        \caption{Solution for single-level SDC (level 2).}
    \end{subfigure}
    \begin{subfigure}[t]{0.48\textwidth}
        \centering
        \includegraphics[width=\textwidth]{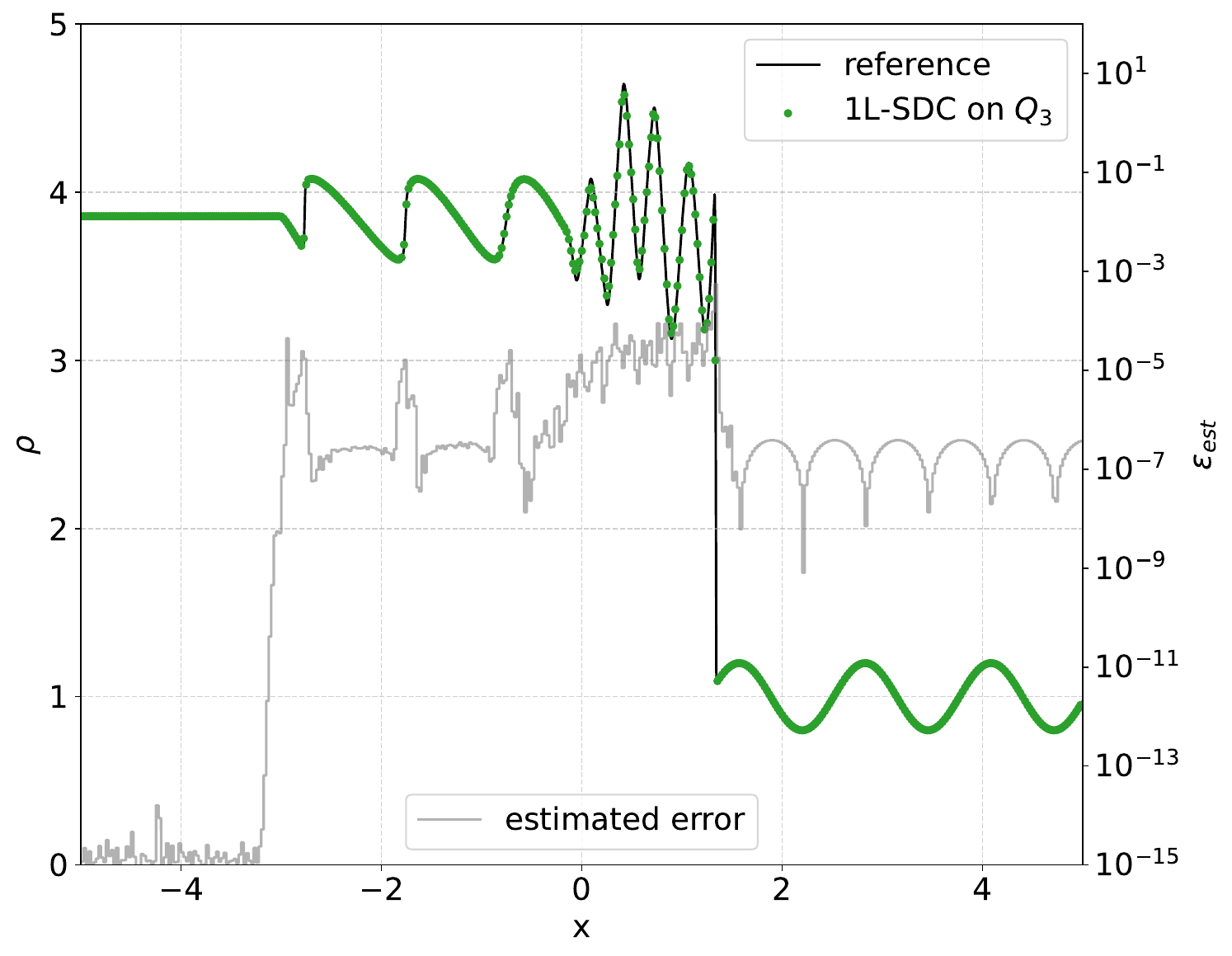}
        \caption{Solution for single-level SDC (level 3).}
    \end{subfigure}
    \hfill
    \begin{subfigure}[t]{0.48\textwidth}
        \centering
        \includegraphics[width=\textwidth]{ch5_shu_osher_t1.500.pdf}
        \caption{Solution for adaptive MLSDC$_{5,7,9}^2$.}
    \end{subfigure}
    \vspace{0.5em}
    \caption{Shock–fluctuation benchmark of Shu \& Osher \cite{TI_Shu1989a}, comparison of single-level SDC solutions with different discretization parameters and space–time adaptive MLSDC$_{5,7,9}^2$-SI(1) solution as well as the estimated error from the novel estimator at $t = 1.5$.
    }
    \label{fig:shu_osher_different_levels}
\end{figure}

The final test case considers the Shu-Osher shock-fluctuation benchmark proposed in \cite[Example~8]{TI_Shu1989a} for the Euler equations with initial conditions
\begin{subequations}
  \begin{alignat}{4}
    &  \rho = 3.857143, &\quad& p = 10.33333, &\quad& v = 2.629369
    && \qquad \text{if} \quad x < -4 \,,\\
    &  \rho = 1 + 0.2\sin 5x, &\quad& p = 1,  &\quad& v = 0
    && \qquad \text{if} \quad x \ge -4 \,.
  \end{alignat}
\end{subequations}
The problem is solved on the domain $\Omega=[-5,5]$ up to $t=1.8$ using a time-step size of $\Delta t=0.001$.
This benchmark involves the interaction of a strong shock wave with a smooth oscillatory density field and is widely regarded as a challenging test case for high-order methods due to the simultaneous presence of discontinuous and smooth solution features.
Since no exact solution is available for this problem, exact-error-based refinement criteria cannot be employed and the assessment must rely on error estimators alone.
Adaptive refinement and coarsening are triggered using thresholds of $5\cdot 10^{-6}$ and $10^{-8}$, respectively.
Throughout this test case, the density $\rho$ serves as the sensor variable for error estimation.
\begin{table}[H]
    \centering
    \begin{tabular}{cccc}
    \toprule
     & \multicolumn{2}{c}{Space} & \multicolumn{1}{c}{Time} \\
    \cmidrule(r){2-3} \cmidrule(r){4-4}
     & $N_{e,l}$ & $P_l$ & $M_l$ \\
    \midrule
    Level 1 & 100 & 8 & 5 \\
    Level 2 & 200 & 8 & 7 \\
    Level 3 & 400 & 8 & 9 \\
    \bottomrule
    \end{tabular}
    \caption{Discretization parameters for the Shu Osher shock tube tests.}
    \label{tab:num_experiments_shu_osher}
\end{table}
\begin{table}[H]
    \centering
    \begin{tabular}{cccc}
    \toprule
     $\Delta t$ & $N_{e}$ & $P$ & $M$ \\
    \midrule
    0.0005 & 400 & 11 & 9 \\
    \bottomrule
    \end{tabular}
    \caption{Discretization parameters for the Shu Osher reference solution.}
    \label{tab:num_experiments_shu_osher_reference}
\end{table}
\begin{figure}[H]
    \centering
    \begin{subfigure}{0.48\textwidth}
        \centering
        \includegraphics[width=\textwidth]{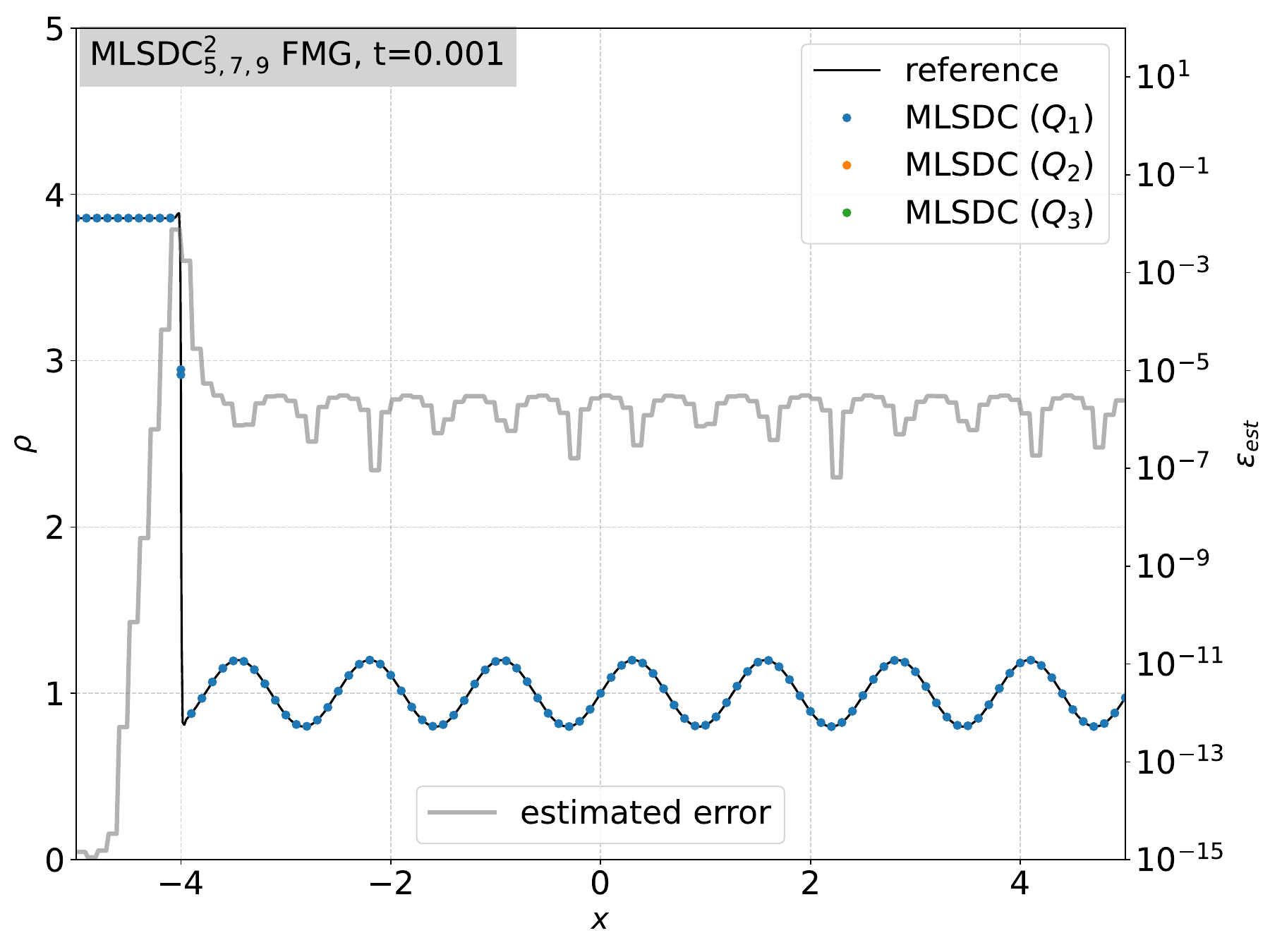}
        \caption{Solution and estimated error at $t = 0.001$.}
    \end{subfigure}
    \hfill
    \begin{subfigure}{0.48\textwidth}
        \centering
        \includegraphics[width=\textwidth]{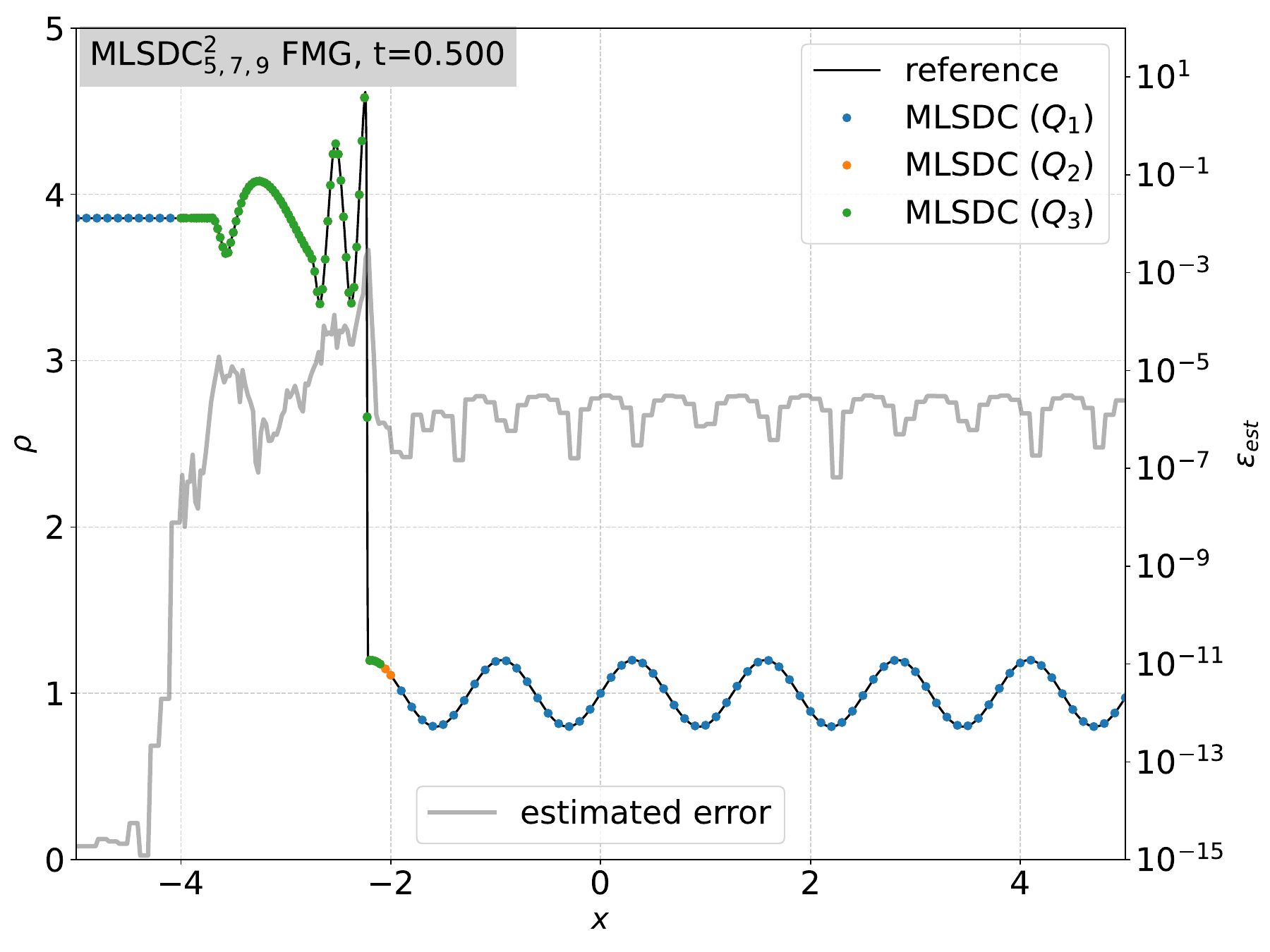}
        \caption{Solution and estimated error at $t = 0.500$.}
    \end{subfigure}
    \begin{subfigure}{0.48\textwidth}
        \centering
        \includegraphics[width=\textwidth]{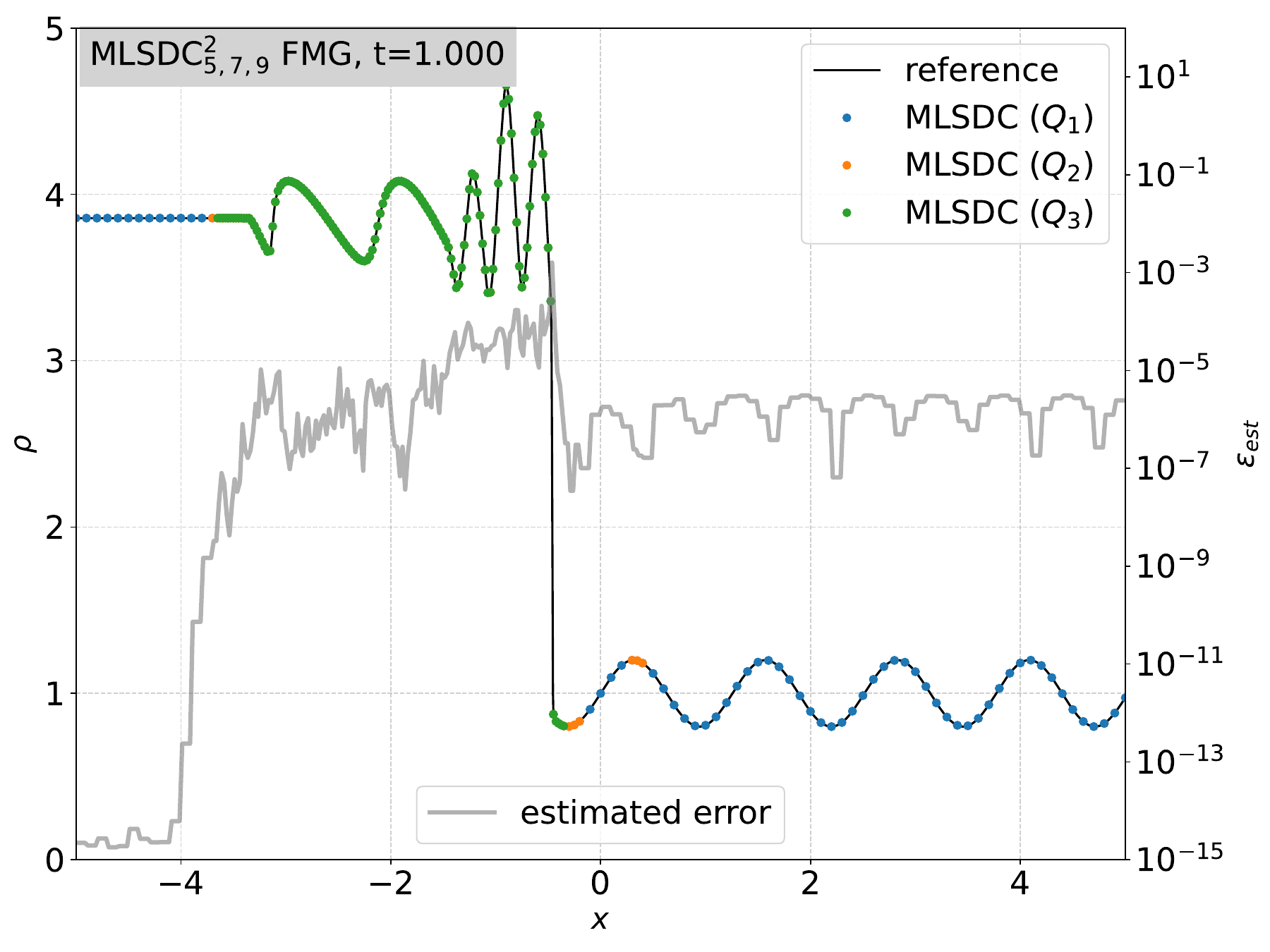}
        \caption{Solution and estimated error at $t = 1.000$.}
    \end{subfigure}
    \hfill
    \begin{subfigure}{0.48\textwidth}
        \centering
        \includegraphics[width=\textwidth]{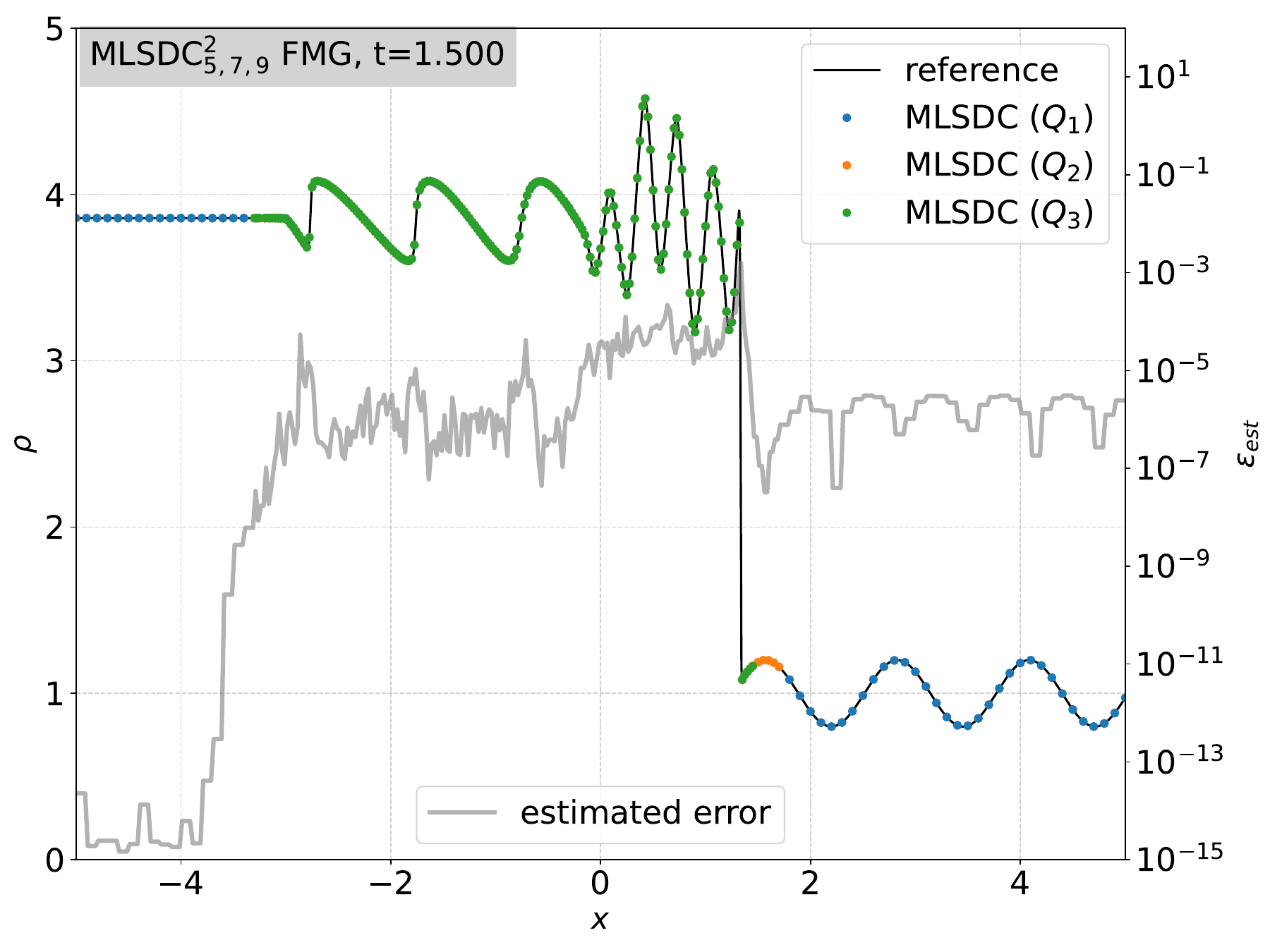}
        \caption{Solution and estimated error at $t = 1.500$.}
    \end{subfigure}
    \medskip 
    \begin{subfigure}{\textwidth}
        \centering
        \includegraphics[width=2.2cm]{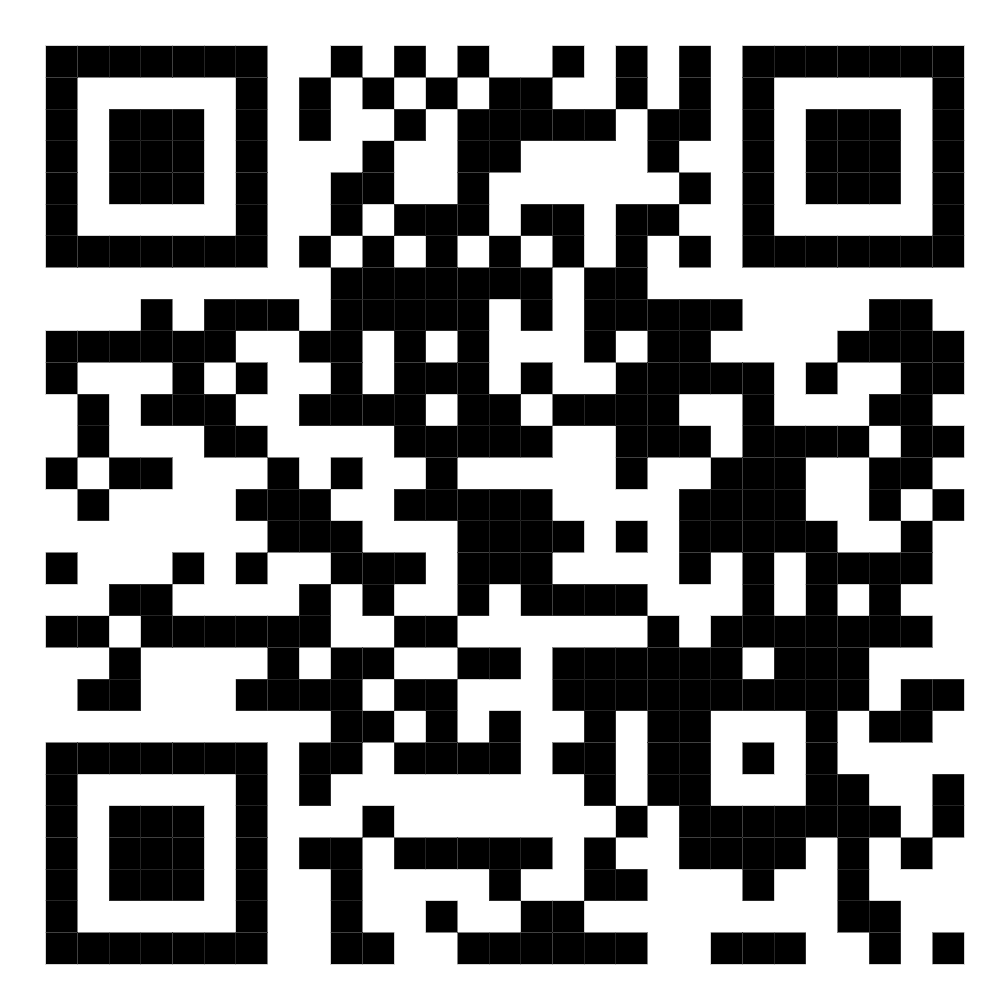}
        \caption{Animation of the full simulation.}
    \end{subfigure}
    \vspace{0.5em}
    \caption{Shock–fluctuation benchmark of Shu \& Osher \cite{TI_Shu1989a} computed with a space–time adaptive MLSDC$_{5,7,9}^2$-SI(1) method as well as the estimated error from the novel estimator for every element. The points indicate element boundaries.
    }
    \label{fig:shu_osher_multilevel}
\end{figure}
Figure~\ref{fig:shu_osher_multilevel} shows the numerical solution for the density $\rho$ at four representative times together with the corresponding element-wise error estimates from the novel estimator. The point markers indicate values located at element boundaries.
The displayed error estimate is obtained from the finest discretization level present within each element. Consequently, the estimated error remains approximately uniform throughout the refined region, which is the intended behavior of the adaptive refinement strategy.
The Henderson estimator was found to be unsuitable for the Shu-Osher benchmark, as its use resulted in numerical instabilities.
The artificial-diffusivity criterion was not considered here, as the results in Subsection~\ref{sec:refine} indicated the criterion's inherent imitations.
As expected, the refinement region follows both the propagating shock and the oscillatory structures of the solution. The largest estimated errors occur in the vicinity of the shock front, indicating that the estimator successfully identifies the most demanding regions of the flow. Since no exact solution is available for this benchmark, a direct comparison with the true error is not possible.
We can moreover observe the adaptive method further coarsens the mesh in regions trailing the front.
To assess the efficiency of the adaptive approach, the adaptive MLSDC$_{5,7,9}^{2}$ method is compared against single-level SDC discretizations corresponding to each level of the multilevel hierarchy. Both solution quality and runtime are considered.
The runtime experiments were performed using the HiSPEET framework compiled with GNU Fortran (GCC) 12.2.0. All simulations were executed in single-threaded mode on an Intel Xeon Platinum 8470 CPU operating at 2.00,GHz. Each configuration was repeated ten times, and the reported runtimes correspond to the average over all runs.
Figure~\ref{fig:shu_osher_different_levels} compares the solutions obtained with the single-level SDC methods of all three levels alone and the adaptive MLSDC method at $t=1.5$. The coarse and intermediate discretizations fail to adequately resolve the fine-scale structures of the solution, whereas the finest single-level discretization produces results comparable to those obtained with adaptive MLSDC.
As shown in Table~\ref{tab:num_experiments_shu_osher_runtime}, the adaptive MLSDC method reduces the runtime by approximately $45\%$ relative to the finest single-level SDC discretization while maintaining a comparable solution quality. These results demonstrate the computational advantages of local space-time adaptivity for this challenging benchmark problem.
The estimated $L^2$ error is included in the figure as well. While the observed differences may appear subtle, please note they are nevertheless significant enough to influence the maximum error attained by the solution.

\begin{figure}[H]
    \centering
    \begin{subfigure}{0.48\textwidth}
        \centering
        \includegraphics[width=\textwidth]{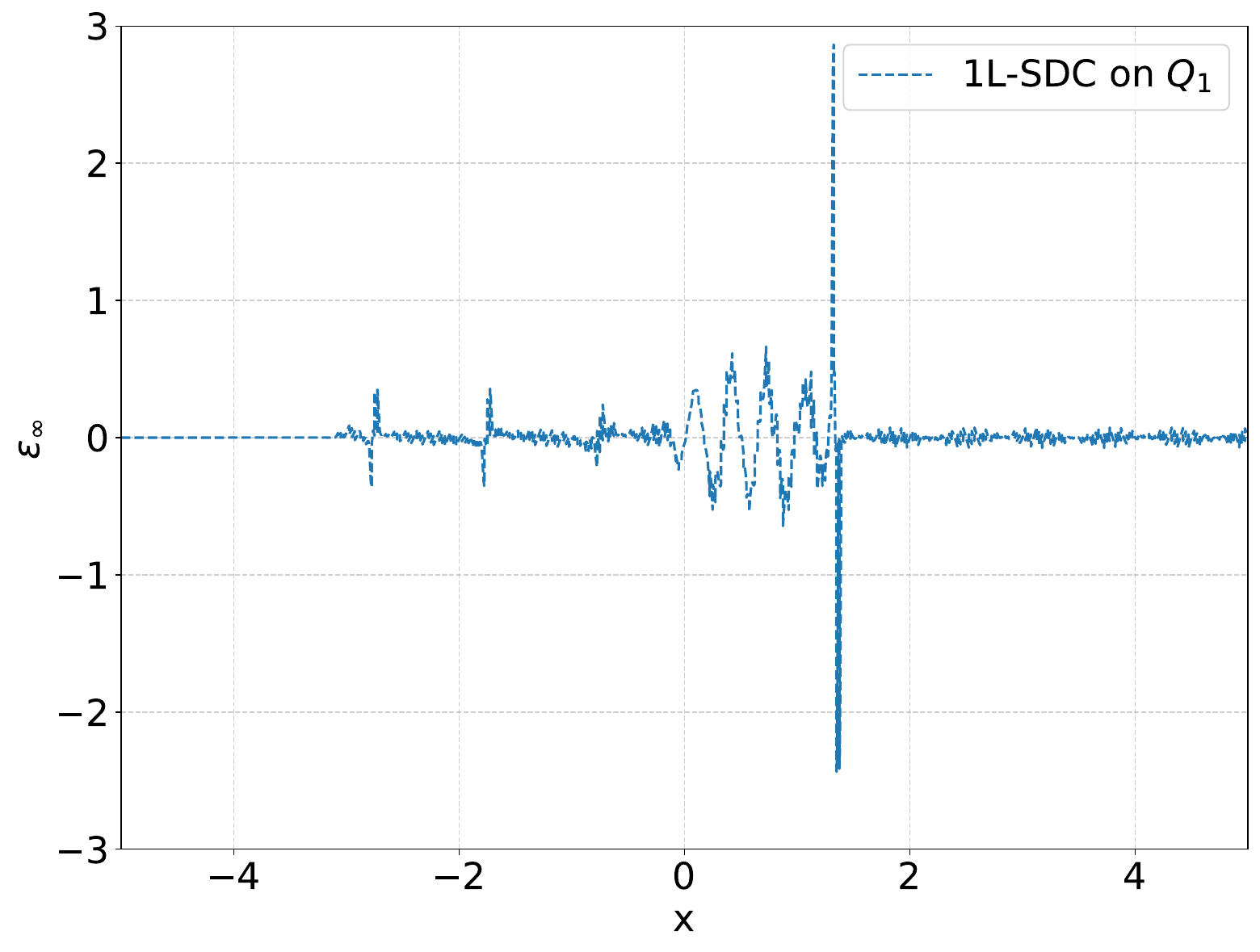}
        \caption{Difference to reference for single-level SDC (level 1).}
    \end{subfigure}
    \hfill
    \begin{subfigure}{0.48\textwidth}
        \centering
        \includegraphics[width=\textwidth]{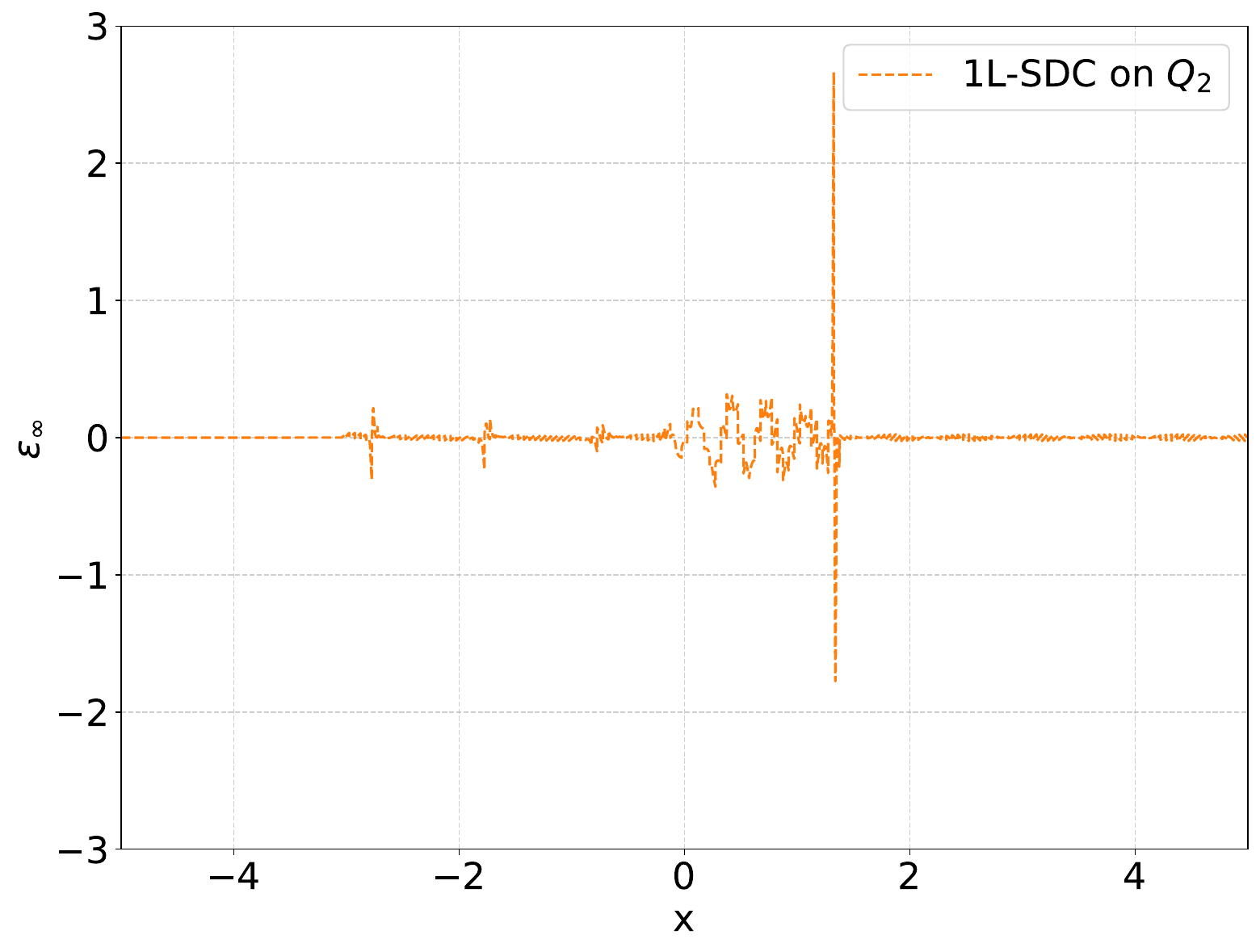}
        \caption{Difference to reference for single-level SDC (level 2).}
    \end{subfigure}
    \begin{subfigure}[t]{0.48\textwidth}
        \centering
        \includegraphics[width=\textwidth]{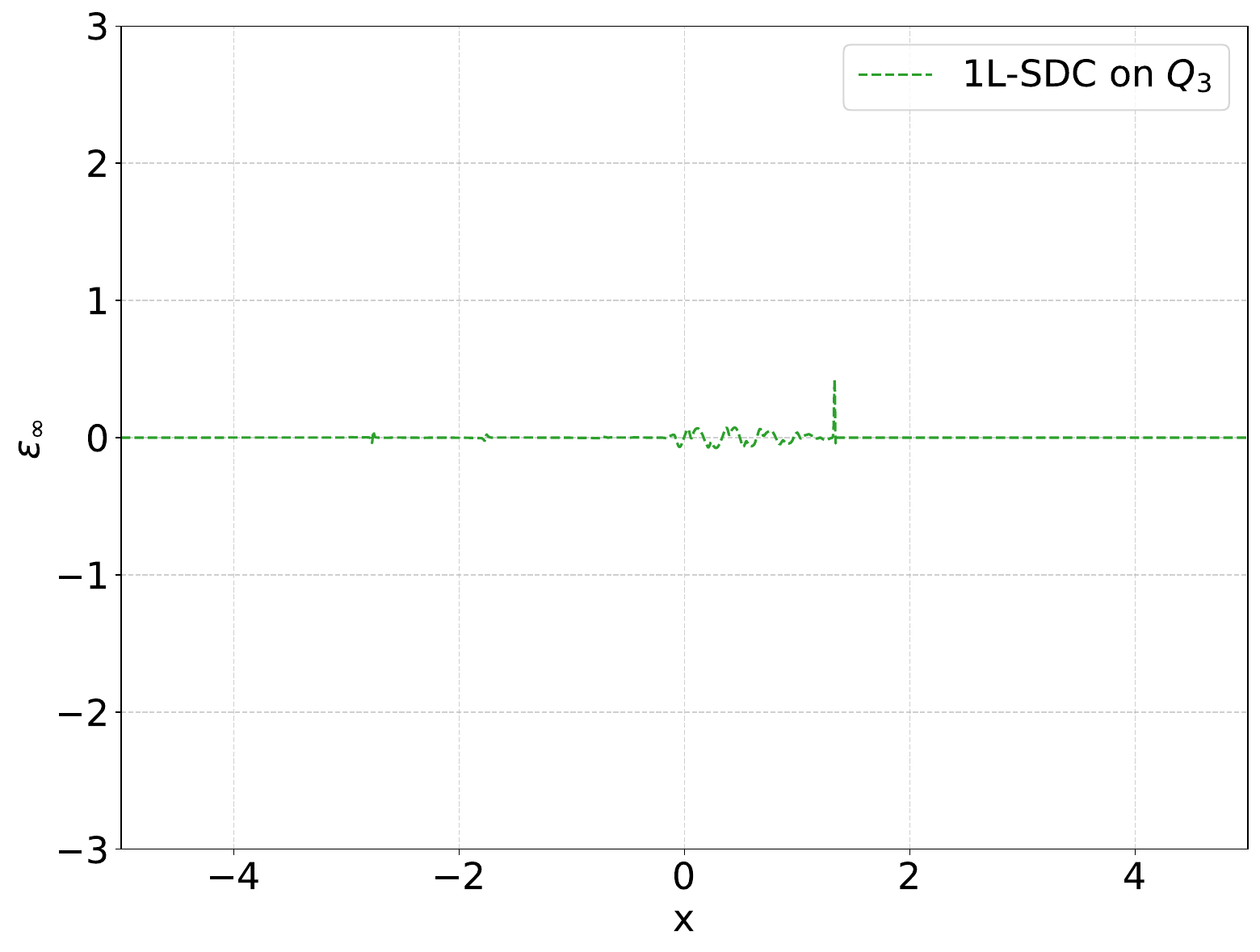}
        \caption{Difference to reference for single-level SDC (level 3).}
    \end{subfigure}
    \hfill
    \begin{subfigure}[t]{0.48\textwidth}
        \centering
        \includegraphics[width=\textwidth]{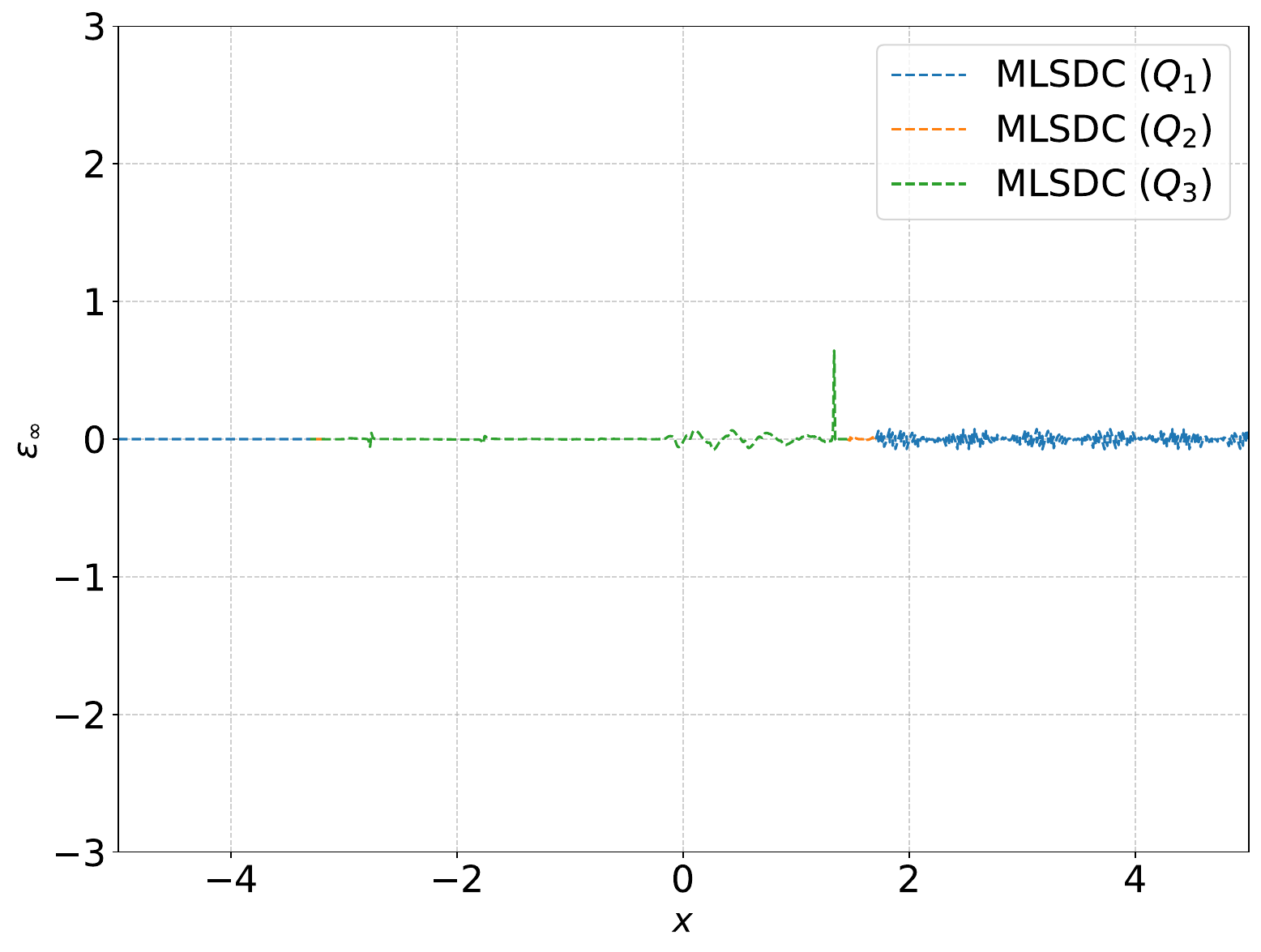}
        \caption{Difference to reference for adaptive MLSDC.}
    \end{subfigure}
    \vspace{0.5em}
    \caption{Shock–fluctuation benchmark of Shu \& Osher \cite{TI_Shu1989a}, comparison of deviations of single-level SDC solutions with different discretization parameters and space–time adaptive MLSDC$_{5,7,9}^2$-SI(1) solution to reference at $t = 1.5$.}
    \label{fig:shu_osher_errors}
\end{figure}

\begin{table}[H]
    \centering
    \begin{tabular}{cccc}
    \toprule
     Method & Runtime [s] & relative Difference \\
    \midrule
    adapt. MLSDC & 3321 & - \\
    SDC (Level 3) & 4819 & +45\% \\
    SDC (Level 2) & 1834 & -45\% \\
    SDC (Level 1) & 677 & -80\% \\
    \bottomrule
    \end{tabular}
    \caption{Runtime comparison for adpative MLSDC$_{5,7,9}^2$ vs SDC on different levels.}
    \label{tab:num_experiments_shu_osher_runtime}
\end{table}
Figure~\ref{fig:shu_osher_errors} depicts the difference between the numerical solutions and a reference solution computed using the discretization parameters listed in Table~\ref{tab:num_experiments_shu_osher_reference} for all four methods under consideration.
The error distributions serves as an addition to the observations made in Figure~\ref{fig:shu_osher_different_levels}. Only the finest single-level SDC discretization and the adaptive MLSDC method achieve sufficiently small errors to accurately resolve the relevant flow structures, whereas the coarser single-level discretizations exhibit noticeably larger deviations from the reference solution.
Although the adaptive strategy is based on a $L^2$ error, we additionally visualize the pointwise deviation from the reference solution. While this representation is not used for a quantitative error assessment, it provides an intuitive view of the error distribution and highlights local differences between the single-level SDC methods and the adaptive MLSDC method. The visualization complements the analysis and helps illustrate the advantages of the adaptive method.


\section{Conclusion}
\label{ch:conclusion}
In this work, a locally space-time adaptive MLSDC method for high-order DG-SEM discretizations has been developed.
The proposed method builds upon existing multilevel SDC approaches, including MLSDC \cite{MG_Speck2014}, MLSDC-SH \cite{MG_Hamon2019a}, AMLSDC \cite{MG_Emmett2019}, and the DG-SEM MLSDC method of \cite{MG_Pfister2025}, which serves as the primary methodological foundation of the present work.
The method extends the MLSDC approach of \cite{MG_Pfister2025} by introducing local space-time adaptivity in order to further improve computational efficiency.
The underlying MLSDC formulation was previously demonstrated to attain convergence orders close to thirteen in both space and time, corresponding to the selected discretization parameters of the study. 
Since both the DG-SEM and MLSDC discretizations are of arbitrary order, these results do not represent an inherent upper limit of the method \cite[Section~8.5]{MG_Pfister2025}.
To this end, adaptive space-time discretizations, the associated adaptive element hierarchy, and the transfer operators as well as their areas of affect required for multilevel computations were introduced and discussed. Based on these ingredients, a locally adaptive FAS collocation formulation and the corresponding MLSDC algorithms were derived.

A key component of the method is the error estimation strategy. For the spatial discretization, a spectral error estimator based on the work of Mavriplis \cite{EE_MavriplisDiss1989} and subsequent developments was adopted.
Since such estimators are not directly applicable to iterative collocation methods, a novel temporal error estimator was developed that accounts for both the superconvergence properties and the iterative nature of SDC methods. Combined with the spatial estimator, it provides a unified space-time error indicator that was shown to yield reliable error predictions across a range of test cases.
To the best of the author's knowledge, the proposed method constitutes the first demonstration of local adaptivity in both space and time within a discretization that simultaneously retains arbitrary-order accuracy in the spatial and temporal dimensions.

The numerical experiments demonstrate the effectiveness of the proposed space-time error estimator. In comparisons with refinement criteria based on Henderson's indicator and artificial diffusivity, the proposed estimator consistently produced competitive or superior results and remained robust for all test cases considered, including the challenging Shu-Osher shock-fluctuation benchmark, for which the alternative criteria proved unsuitable.
For the Burgers moving-front problem, substantial reductions in computational cost were observed even in one spatial dimension. The globally refined MLSDC method already outperformed the IMEX-RK ARS(4,4,3) scheme over most of the investigated error range (compare also \cite[Section~8.4]{MG_Pfister2025}), while the adaptive MLSDC method provided additional efficiency gains through localized refinement.
The parameter studies further indicate that the performance of the adaptive MLSDC method depends strongly on the refinement and coarsening thresholds, which directly determine the extent of the refined region. 
Consequently, the systematic optimization of refinement criteria and threshold selection represents a promising direction for future research.

The proposed  adaptive MLSDC method was shown to robustly resolve the challenging Shu-Osher shock-fluctuation benchmark with the novel estimator, demonstrating its applicability to complex flow problems involving both discontinuities and smooth structures.
For this test case, the adaptive three-level MLSDC method attained the accuracy of the finest-level discretization while reducing the computational cost by approximately $45\%$.
Given that the computational savings originate from localized space-time refinement, the potential efficiency gains are expected to increase considerably when extending the method from one-dimensional to three-dimensional applications.

Several directions for future research remain.
The extension of the proposed adaptive MLSDC method to three-dimensional problems constitutes a particularly promising direction for future work. Owing to the localized nature of the adaptive refinement strategy, the computational savings are expected to increase substantially with problem dimensionality.
In the present work, refinement is performed in a coordinated space-time manner, with the discretization parameters on the finer levels prescribed a priori. This choice was made deliberately in order to isolate and investigate the effects of locally increased space-time resolution.
However, the proposed methodology is not restricted to coupled space-time refinement. Since the spatial and temporal error contributions are estimated separately before being combined.
Future extensions could exploit this information to selectively introduce additional levels for either spatial or temporal refinement, depending on the dominant source of error. In combination with the favorable stability properties of the employed semi-implicit integrators, such anisotropic space-time adaptation represents a promising direction for future research.
Although the proposed adaptive MLSDC method is technically designed to support temporal $h$-adaptivity, a sufficiently robust implementation has not yet been achieved. The development and assessment of locally adaptive temporal mesh $h$-refinement strategies therefore remain subjects of ongoing investigation.
Another promising direction is the construction of MLSDC schemes based on alternative Runge-Kutta collocation formulations with improved stability properties (TVD or SSP), which may further enhance the robustness of the method for strongly convection-dominated problems.


\section*{Data Availability}

The datasets generated during the current study are available from the corresponding author on reasonable request.

\section*{Funding}

Funding by German Research Foundation (DFG) in frame of the project STI 157/9-1 is gratefully acknowlged. The authors would like to thank ZIH, TU Dresden, for the provided computational resources. Open Access funding is enabled and organized by Projekt DEAL.

\section*{Ethics declarations}

The authors declare that they has no known competing financial interests or personal relationships that could have appeared to influence the work reported in this paper.


\bibliographystyle{abbrvnat}
\bibliography{BM,FE,FV,SE,MG,TI,EE,ST}


\end{document}